\documentclass[12pt]{article}
\usepackage{amsmath}
\usepackage{amsthm}
\usepackage{amscd}
\usepackage{amssymb,amsfonts}
\usepackage[backend=biber,style=alphabetic,maxnames=5]{biblatex}
\usepackage{listings}

\usepackage{hyperref}
\usepackage{latexsym}
\usepackage{todonotes}
\usepackage{nicefrac}
\usepackage{epsfig}
\usepackage{stmaryrd}
\usepackage{setspace}
\usepackage{enumerate}
\usepackage[all]{xypic}
\usepackage{bbm,ifpdf,tikz}
\usepackage{verbatim}
\ifpdf
\usepackage{pdfsync}
\fi
\usepackage[T1]{fontenc}
\usepackage[utf8]{inputenc}
\usepackage{float}

\newtheorem{theorem}{Theorem}[section]

\newtheorem{conjecture}[theorem]{Conjecture}

{
\theoremstyle{definition}
\newtheorem{definition}[theorem]{Definition}
\newtheorem{example}[theorem]{Example}

\newtheorem{remark}[theorem]{Remark}
\newtheorem{openthm}{Open Question}
}

\newcommand{\excise}[1]{}

\renewcommand{\and}{\qquad\text{and}\qquad}

\newcommand{\prufer}{Pr\"ufer } 
\newcommand{\orderedges}{\texttt{order\_edges }}

\title{An AI enhanced approach to the tree unimodality conjecture}
\author{Eric Ramos and Sunny Sun}

\begin{document}
\maketitle

\begin{abstract}
Given a graph $G$, its independence sequence is the integral sequence $a_1,a_2,\ldots,a_n$, where $a_i$ is the number of independent sets of vertices of size i. In the late 80's Alavi, Erd\"os, Malde, Schwenk showed that this sequence need not be unimodal for general graphs, but conjectured that it is always unimodal whenever G is a tree. This conjecture was then naturally generalized to claim that the independence sequence of trees should be log concave, in the sense that $a_i^2$ is always above $a_{i-1}a_{i+1}$. This conjecture stood for many years, until in 2023, Kadrawi, Levit, Yosef, and Mizrachi proved that there were exactly two trees on 26 vertices whose independence sequence was not log concave. In this paper, we use the AI architecture PatternBoost, developed by Charton, Ellenberg, Wagner, and Williamson to train a machine to find counter-examples to the log-concavity conjecture. We will discuss the successes of this approach -- finding tens of thousands of new counter-examples to log-concavity with vertex set sizes varying from 27 to 101 -- and some of its fascinating failures.
\end{abstract}

\section{Introduction}
Given a graph $G$, its independence sequence is the integral sequence $a_1,a_2,\ldots,a_n$, where $a_i$ is the number of independent sets of vertices of size $i$. In the late 80’s Alavi, Erd\"os, Malde, Schwenk showed that this sequence need not be unimodal for general graphs, but conjectured that it is always unimodal whenever G is a tree\cite{Alavi1987VertexIndependence}. 

In the last twenty years or so, a stronger condition, log-concavity -- stipulating that $a_i^2$ is always above $a_{i-1}a_{i+1}$ -- was considered in an attempt to prove the tree unimodality conjecture. Though this condition is strictly stronger than log-concavity, and thus ostensibly harder to prove, it is also much more amenable to geometric arguments, thereby making it the first plan of attack in many modern treatments \cite{huh2022combinatorics}. In fact, this stronger version of the conjecture was shown to hold for all trees of at most 25 vertices \cite{yosef2022unimodality}. In 2023, however, using improved computational power and a considerably more efficient algorithm, Kadrawi, Levit, Yosef, and Mizrachi proved that there were exactly two trees on 26 vertices whose independence sequence was not log concave. They also showed how these two examples could be generalized to create a few families of trees whose members are all not log concave. Finally, in early 2025, Galvin provided a family of trees with the property that for any chosen positive integer $k$, there is a member T of the family where log-concavity breaks at index $\alpha(T) - k$, where $\alpha(T)$ is the independence number of T. 

Motivated by the the recent explosion in the use of machine learning techniques to augment mathematical discovery \cite{davies2021advancing,shehper2024makes,romera2024mathematical,ellenberg2025generative} and the userfriendlyness of the machine learning architecture Patternboost \cite{charton2024patternboost}, developed by Charton, Ellenberg, Wagner, and Williamson, we trained a machine to cook up a range of counterexamples to log-concavity in the hope that doing so may illuminate (in the positive or negative) the original tree unimodality conjecture. In total, our methods have thus far uncovered tens of thousands of new counter-examples to the log-concavity conjecture. We explicitly present a number of these counter-examples in the Appendix. In our github repo -- found at the url \url{https://github.com/ericgramos/TreeUnimodalityPatternBoost} -- one will find a folder with the complete output of one experiment we ran on trees with 60 vertices. Within that folder, there is a text file containing around 35,000 trees on 60 vertices whose independence sequence is not log concave.

Before this work, all of the handful of known families of counter-examples looked very similar to one another (see Figure \ref{fig: first counters} for the general ``flavor" we are talking about). One important consequence of the present work is that it doesn't seem as though the similarities in these early counter-examples were indicative of general trends. Indeed, as one allows the number of vertices of the tree to grow, our methods were able to produce trees that, to our eyes, seem to be essentially generic, albeit perhaps with an excess of vertices of degree 2. The Appendix of this work illustrates just over a dozen of the counter-examples that we discovered across a variety of vertex sizes. We hope that our curated list gives one a sense of how variable these counter-examples can be.

This paper is organized as follows: we begin by detailing the history of the tree unimodality conjecture, and present background related with the PatternBoost architecture. Following this, we spend the bulk of the paper explaining all of our choices when it came to building the machines we used to find all of the new counter-examples. To finish the paper, we review all of the various results and failures of our approach.

\subsection{Directions for future study}

There are a number of future directions that one can take based on the results of this one.

\begin{itemize}

\item All of the counter-examples we were able to find using our methods exhibit the breakage of log-concavity within 3 indices of the top of the independence sequence. We touch upon why this might be in Sections \ref{sec:alpha} and \ref{sec:N/2-1}, though this is largely speculative. Are there modifications one can make to our machine to find counter-examples whose breakage happens at different indices? For instance, we were not able to find counter-examples with breakage three indices below the independence number until we began punishing the path (see  Example \ref{ex:alpha-3})

\item Related with the above point, we were also never able to find counter-examples where log-concavity breaks more than once. The question of whether it is even possible for log-concavity to break twice in this context was originally posed by Galvin in \cite{galvin2025trees}. The reward function we used to train our models targets individual indices for breakage. Is it possible that a more well rounded scoring function can be used to force breakage in multiple spots?

\item Using a reward function that teaches the machine to modify the independence sequence more ``globally" rather than ``locally" might be useful to attack the more general unimodality conjecture. It remains unclear at this point what such a reward function would look like.

\end{itemize}

\subsection{Standardizing notation}

Here we collect certain notation we will use throughout the work.

\begin{itemize}
    \item $N$ -- the number of vertices of the tree being discussed. Equivalently, two more than the length of the \prufer code of the tree.
    \item $N/2$ -- this \textbf{always} refers to the floor of $N/2$ when $N$ is odd.
    \item $\alpha, \alpha(G),\alpha(T)$ -- the independence number of the graph or tree being discussed.
    \item $a_i$ -- the number of independent sets of the graph under discussion with size $i$.
    \item $\mathbf{3_{44}}$,$\mathbf{3^*_{34}}$ -- the original two counter-examples to log-concavity on 26 vertices. See \ref{fig: first counters} for illustrations of both.
    \item $I_G(x)$ -- the independence polynomial of graph $G$, $I_G(x) = \sum_{i=0}^\alpha a_ix^i$
\end{itemize}

The code accompanying this paper can be found at: \\ \url{https://github.com/ericgramos/TreeUnimodalityPatternBoost}.
The original PatternBoost code can be found here:
\url{https://github.com/zawagner22/transformers_math_experiments.}

\section*{Acknowledgements}
This project was completed as part of the NYC Discrete Math REU. We thank Adam Sheffer for running this program and Baruch College for their facilities. The necessary super computing power for our experiments was provided by a collaboration between the Stevens Institute for Artificial Intelligence (SIAI), and the Stevens Laboratory Artificial Intelligence in Mathematics Education. The authors would like to send thanks to Harvard CSMA's Math and Machine Learning reunion workshop, where the results of this work were originally announced.

Both authors would like to send their thanks to Francois Charton, Jordan Ellenberg, Steven Heilman, and Greta Panova for numerous thought provoking conversations. Big thanks must also be sent to David Galvin and Ferenc Bencs, who reached out to us after the first draft of this manuscript was posted, and let us know about their family of examples for which log concavity can be made to break multiple times.

The first author is especially grateful to Professor Ellenberg, who was kind enough to never stop responding to his emails. He is also thankful for Wanjia Guo, who taught him how to use command line. The first author was partially supported through NSF grant DMS-2452031.

The second author is grateful to Shiyue Li for pointing her to the paper \cite{galvin2025trees} that was inspirational to this project. She would also like to thank her peers in the REU for sharing not only their insights but also moments of encouragement, laughter, and support that made this experience truly meaningful.

\section{Background}

\subsection{The tree unimodality conjecture and log-concavity}

In this section we will run through the various background necessary for the tree unimodality conjecture.

\begin{definition}
Given any simple graph $G=(V,E)$, an \textbf{independent set} $\{v_1,\dots,v_n\}$ is a set of pairwise non-adjacent vertices. In other words, $(v_i,v_j) \notin E$ for all $i,j$.

The \textbf{independence number $\alpha(G)$} for a graph $G$ is the cardinality of its largest independent set. 

Let $a_i$ denote the number of independent sets of cardinality $i$ in a given graph $G$. The \textbf{independence sequence} of a graph $G$ is given by $\{{a_i}\}^{\alpha(G)}_{i=0}$ where $a_0 = 1$. We frequently think of the $a_i$ as coefficients of a polynomial, which we denote $I_G(x)$, called the \textbf{independence polynomial} of the graph $G$.
\end{definition}

It is well known that computing independent sets of a general graph is NP hard. Despite this, there are a number of families for which the problem turns out to be considerably more tractable. One such family, and the one most relevant to us, are trees. This largely follows from the natural recurrence,
\[
I_G(x) = I_{G-v}(x) + x\cdot I_{G-N(v)}(x),
\]
where $v$ is any vertex and $N(v)$ is its (closed) neighborhood. In the case where $G$ is a tree, or forest, the lack of cycles means that this recurrence can be performed dynamically, minimizing the necessary computations to be linear in the number of vertices. This algorithm, which was pivotal in our own computations, was explored in \cite{Kadrawi2023OnComputing}.

As is common in the study of finite combinatorial motivated sequences, it is natural to ask whether the independence sequence is unimodal. Indeed, it was shown by Hamidoune \cite{hamidoune1990numbers} that the independence sequence of any claw-free graph was always unimodal. Generally speaking, however, it was proven in \cite{Alavi1987VertexIndependence} that the independence sequence of graphs could be arbitrarily far from unimodal. In that same work, it was conjectured that the sequence is unimodal in the case of trees.

\begin{conjecture}[The Tree Unimodality Conjecture, \cite{Alavi1987VertexIndependence}]For any tree, the sequence of independence numbers \( a_0, a_1, \ldots, a_\alpha \) is unimodal: there exists an index \( m \) such that
  \[
  a_0 \leq a_1 \leq \cdots \leq a_m \geq a_{m+1} \geq \cdots \geq a_\alpha.
  \]
\end{conjecture}

In view of Hamidoune's theorem, this conjecture suggests that unimodality of the independence sequence will hold at the two extremes of claw-free and ``claw-full," while the results of \cite{Alavi1987VertexIndependence} imply that things are much more complicated in the interim.

Recent advancements in the geometry and combinatorics of matroids (see \cite{huh2022combinatorics} for a survey) have shone a light on a property that is strictly stronger than unimodality of a sequence (with no internal zeros), log-concavity.

\begin{definition}
A sequence of nonnegative real numbers $a_0, a_1, \ldots, a_n$  is called \textbf{log-concave} if
\[
a_i^2 \geq a_{i-1} a_{i+1}
\quad \text{for } 0 < i < n.
\]
\end{definition}

In \cite{chudnovsky2007roots}, Chunovsky and Seymour expanded Hamidoune's result to show that the independence sequence of claw-free graphs was log-concave. More recently, Li was able to prove an equivariant version of this statement as well \cite{LiEquivariant}. Somewhat optimistically for the tree unimodality conjecture, using an algorithm based on a database of nonisomorphic unlabeled trees due to Yosef, Mizrachi, and Kadrawi \cite{yosef2022unimodality}, Redcliffe verified the log-concavity of tree independence sequences with up to 25 vertices \cite{ball2022independent}.

\begin{figure}
    \centering
    \includegraphics[width=0.5\linewidth]{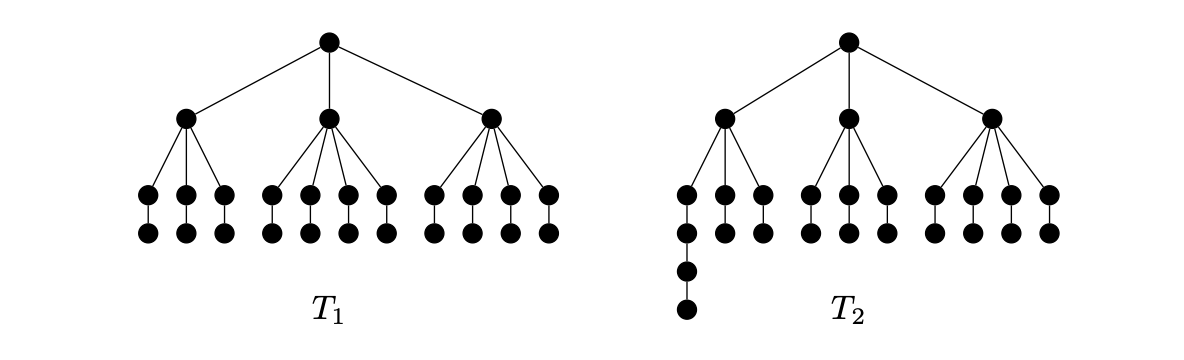}
    \caption{The two counter-examples, $T_1 = \mathbf{3_{44}}$ and $T_2 =\mathbf{3^{*}_{34}}$ to log-concavity with 26 vertices. These images have been taken from their source \cite{Kadrawi2023OnComputing}}
    \label{fig: first counters}
\end{figure}

In 2023, however, Kadrawi, Levit, Yosef, and Mizrachi used a dynamic programming approach to extend these calculations and concluded that there were exactly two trees on 26 vertices whose independence sequence was not log-concave \cite{Kadrawi2023OnComputing}. These trees can be found in Figure \ref{fig: first counters}. So it is convenient to refer back to these in future discussions, we follow the naming conventions of \cite{Kadrawi2023OnComputing}, and refer to the left tree $T_1$ as $\mathbf{3_{44}}$, and the right tree $T_2$ as $\mathbf{3^{*}_{34}}$.

In followup work, Kadrawi, Levit, Yosef, and Mizrachi were able to generalize these two counterexamples to seven families of trees with non-log-concave independence sequences \cite{kadrawi2023independence}. It was noted in that work that all members of these families experienced the breakage of log-concavity at only one point in the sequence, and that this point was only ever in position either one or two below the top of the sequence.

This then leads to the question of \emph{where} in the independence sequence can breakage of log-concavity happen. Specifically, Kadrawi and Levit conjectured that for any arbitrary $k$, there should exist a tree $T$ for which log-concavity may be broken at the $\alpha(T)-k$ position \cite{kadrawi2023independence}.

Galvin resolved the conjecture through the construction of a family of trees where for any integer $k \geq 1$, there is a tree $T(k)$ in the family for which log-concavity of its independence sequence breaks at $\alpha(T)-k$ \cite{galvin2025trees}. Galvin's work has kept two notable important doors open, however. Firstly, in all counter-examples thusfar found, including \cite{Kadrawi2023OnComputing, kadrawi2023independence,galvin2025trees}, and even the present work, log-concavity breakage only seems to happen at most once in the independence sequence.

\begin{remark}
    Immediately after the public release of this work, David Galvin and Ferenc Bencs reached out to the authors to inform us they had found a family of counter-examples where log concavity breaks in multiple indices. Interestingly, even the smallest examples that they constructed had well over 100 vertices, an amount that is out of the range of the current work. The conclusions of both this work and their work strongly suggests that if something like unimodality does eventually fail, it would likely happen at a vertex count that is far beyond our computational abilities. 
\end{remark}

Secondly, in all of the counter-examples of prior works, the trees had a very similar ``flavor." These trees are constructed by starting with a root vertex, which is given some number of children. Each of these children is then given a uniform number of children of their own, with the added twist that these grandchildren of the root are paths rather than vertices. One of the primary contributions of this work is to illustrate that this flavor seems to be a consequence of the ``law of small numbers" rather than anything deeply related with the structural properties of trees whose independence sequence is not log concave. One can look to any number of the examples in the Appendix of this work, or the tens of thousands of examples on 60 vertices present in the github repo of this work, to see that counter-examples can come in a wide variety of forms. Indeed, the only consistent theme that is noticeable at this point across all counter-examples is that they tend to have a significant number of vertices of degree 2, especially vertices of degree 2 that are connected to leaves of the tree.

\begin{openthm}
    What are the structural properties of trees whose independence sequence fails to be log concave? In the large vertex limit, is the proportion of trees that fail to be log concave positive?
\end{openthm}

\subsection{PatternBoost}

We used the machine learning architecture PatternBoost developed by Charton, Ellenberg, Wagner, and Williamson. PatternBoost is an accessible tool with successful applications to several problems in extremal combinatorics \cite{charton2024patternboost}. We refer the reader to the original paper for details on the machine learning aspects of how PatternBoost functions. For our purposes, we will only give a general idea about the philosophies underlying PatternBoost

The most key feature of PatternBoost is its local-global iteration, which consists of a local phase with a simple, often times greedy, search algorithm and a global phase of training a transformer. After creating a database and defining a \textbf{scoring function}, which evaluates how good a construction is, PatternBoost starts its \textbf{local search}: it takes graphs in the database, makes small local modifications with the goal of improving scores, and outputs graphs with scores not lower than that of the input graph.

The global phase involves training a transformer on the best constructions discovered during the local phase. PatternBoost uses the transformer implementation Makemore, developed by Andrej Karpathy for its accessibility and flexibility \cite{karpathy2022makemore}. In each generation, the transformer is trained with the provided database, and is then asked to produce similar constructions. We filter out the invalid constructions, and plug the rest to the local search algorithm for local optimization. Finally, we select the graphs with best scores in the database as the training set for the next generation and continue the process with a finetuned database. Each generation of local to global is known as an \textbf{epoch} of the process.

In the following sections we will detail how we chose to represent trees for the purposes of the machine, the scoring function meant to capture the failing of log-concavity, and the local search algorithm that was implemented to improve this score. We will also go into much more details in relation to the parameters of our experiments in Section \ref{sec:parameters}.

\section{Building the machine}

\subsection{The choice of representation and scoring function}

Aware that the transformer in PatternBoost takes one-dimensional inputs, we chose \prufer codes to encode trees. Note that while PatternBoost comes equipped with a tokenizer which can in principle convert any input into a short string that the transformer would understand, using \prufer codes has the added benefit that this phase of the process is one-to-one. There is no risk of information loss because of the need to combine larger chunks of information into a single token.

\begin{definition}
    A \textbf{\prufer code} is a string of length $n-2$ on the alphabet $\{1,2,\dots,n\}$. By a famous construction of \prufer \cite{Prufer1918}, these codes are in natural bijection with spanning trees of the complete graph $K_n$. In other words, a \prufer code can be thought of as a tree on $n$-vertices whose every vertex has a unique label in the set $\{1,\ldots,n\}$, up to label preserving isomorphism. We will refer to the code and the tree corresponding to this code interchangeably going forward. 
\end{definition}

\prufer codes offer a compact representation of trees, as their length grows linearly with the number of vertices. It also offers a few other conveniences in implementing PatternBoost for our specific project: To verify whether a construction is a valid tree, we just need to check the length of its \prufer code and if the indexes are out of bounds. In addition, a path graph (i.e. a graph whose vertices can be listed in order $v_1, v_2, \dots, v_n$ such that the edges are $(v_i,v_{i+1})$ where $i=1,2,\dots,n-1$) corresponds to \prufer codes with no repeats, whereas a star corresponds to a \prufer code using only one member of the alphabet. This allows us to identify paths and stars with virtually zero computational baggage, which will be convenient in our local search (described in the next section).

We believe that it is also the case that having a more direct representation, with tokens directly corresponding to labeled vertices in trees, could have also increased the learning speed and learning outcomes for the machine.

Of course, with all of the above said, there are certainly possible down-sides to this choice of representation. For instance, while it is known that the number of isomorphism classes of trees grows roughly exponentially in the number of vertices -- no faster than $4^n$ -- the number of \prufer codes grows superexponentially -- precisely $n^{n-2}$. Moreover, because the labels on the vertices contribute nothing to the computation of the independence polynomial, all of these extra trees are essentially just noise. Moreover, once again because the labels play no role in the computation of the independence polynomial, there is some risk in the machine realizing that once it obtains a high enough score, it can repeatedly reach this score by just permuting the labels.

While both of these concerns are legitimate, we ultimately found that neither prevented us from using PatternBoost to create counter-examples to log-concavity. In fact, as we will discuss later, it was very rare for the machine to find the same tree with non log concave independence sequence more than once. It never learned to ``cheat" by permuting labels in that way.

Moving on to our scoring function, we decided to define our scoring function as the literal measure of the failure of log-concavity at a chosen fixed index $i$. To be more precise, given the input index $i$,
\[
\text{score}_i = a_{i-1} \cdot a_{i+1}-a_i^2
\]
where, as in the background section, $a_i$ is the number of independent sets with size $i$. Observe that, by definition of log-concavity, any given input tree's independence sequence will have breakage of log-concavity in its $i$-th index if and only if the above score is positive.

\begin{remark}
    As a side note, we were surprised to see just how varied the values of this scoring function could be. For instance, looking at the examples in the Appendix, one will find trees with scores ranging from around $10^{14}$, all the way down to the smallest possible positive score of 1.
\end{remark}

Note that this score function does not look at the entire independence sequence, but rather only a single index. We believe that this narrow focus allowed the machine to more easily find trees with log-concavity breakage at the individual indices requested. However, we must also acknowledge that this may have prevented the machine from ever truly learning how to manipulate the more global ``shape" of the independence sequence. In particular, in all of our experiments, the machine never managed to produce an example where log-concavity broke twice or more. Of course, this may also be because breaking twice is impossible!

Further note that through our experimentation we have qualitative evidence that our choice of score function limits the scope of the transformer's perspective. In one of our early experiments, we ran our methods on trees with 101 vertices. We did this not realizing that at this vertex count the sizes of the the integers $a_i$ could be higher than Julia's maximum integer size of $2^{64}$. In particular, the independence polynomial computations being done were producing total non-sense around the middle of the polynomial. Despite this, the counter-examples that the machine produced continued to be counter-examples when the overflow was corrected. To be clear, the index of the sequence where we were looking for breakage never overloaded, but other indices did. This suggests to us that the machine really paid no mind whatsoever to the indices of the independence polynomial outside of the three it needed to compute the scoring function.

The definition of our scoring function puts the burden of deciding which index to look at on the user. For any fixed index $i$, independent of the number of vertices $N$, the value of $a_i$ will, on average, look almost identical to $\binom{N}{i}$ when $N$ is big enough. In particular, one does not expect for indices independent of $N$ to produce counter-examples once $N$ gets big enough. Therefore, it seems prudent to choose an index that varies in some way with $N$. As mentioned in the background sections, the original counter-examples of \cite{Kadrawi2023OnComputing,kadrawi2023independence} all had breakage in either index $\alpha(G)-1$ or $\alpha(G)-2$. Because of this, one may think to consider scores in indices of the form $\alpha(G)-k$,where $k$ is a fixed pre-determined integer. 

However, in experiments in which we tried to actively manipulate a quantity relative to the independence number during optimization, we frequently found the machine getting lost. See Remark \ref{rmk:smallBanFail} for more on this. Suffice to say that the machine performed considerably worse whenever we chose an index defined relate to the entered tree's independence number.

Instead of this, we noticed that in the counter-examples of \cite{Kadrawi2023OnComputing,kadrawi2023independence}, the breakage of log-concavity always happened at the index equal to (the floor of) half the number of vertices of the tree\footnote{Throughout this work, whenever we discuss $N/2$ as being an index, we are implicitly talking about the floor of $N/2$}. We therefore spent many of our experiments looking at that particular index. In fact, with our methods we were ultimately only able to find counter-examples with breakage at the indices of $N/2$ and $N/2-1$. We will discuss this more in our final section.

Note that, using the work of Galvin \cite{galvin2025trees}, for instance, it is known that log-concavity breakage should exist within the vertex counts that we checked at indices other than $N/2$ and $N/2-1$. Despite this, we were never able to find any of them. In fact, even finding counter-examples at the index $N/2-1$ required specific engineering. What this all points to is the fact that for some reason these methods had a relatively easy time manipulating independent sets at size around half the number of vertices, and a much harder time otherwise. We will explain why we think this is happening in later sections, once we establish more about how we locally optimized given trees.

\subsection{The local search algorithm}\label{sec:local}

Throughout its applications to several problems in extremal combinatorics, PatternBoost has performed well given very simple local search algorithms with is local-global iterations \cite{charton2024patternboost}. In keeping with the philosophy of Sutton's now influential article, ``The Bitter Lesson," it frequently matters more to these AI assisted processes that things be taken to scale, rather than things be done in the most perfectly optimal or sophisticated way \cite{sutton2019bitter}.

\begin{figure}
    \centering
    \includegraphics[width=0.5\linewidth]{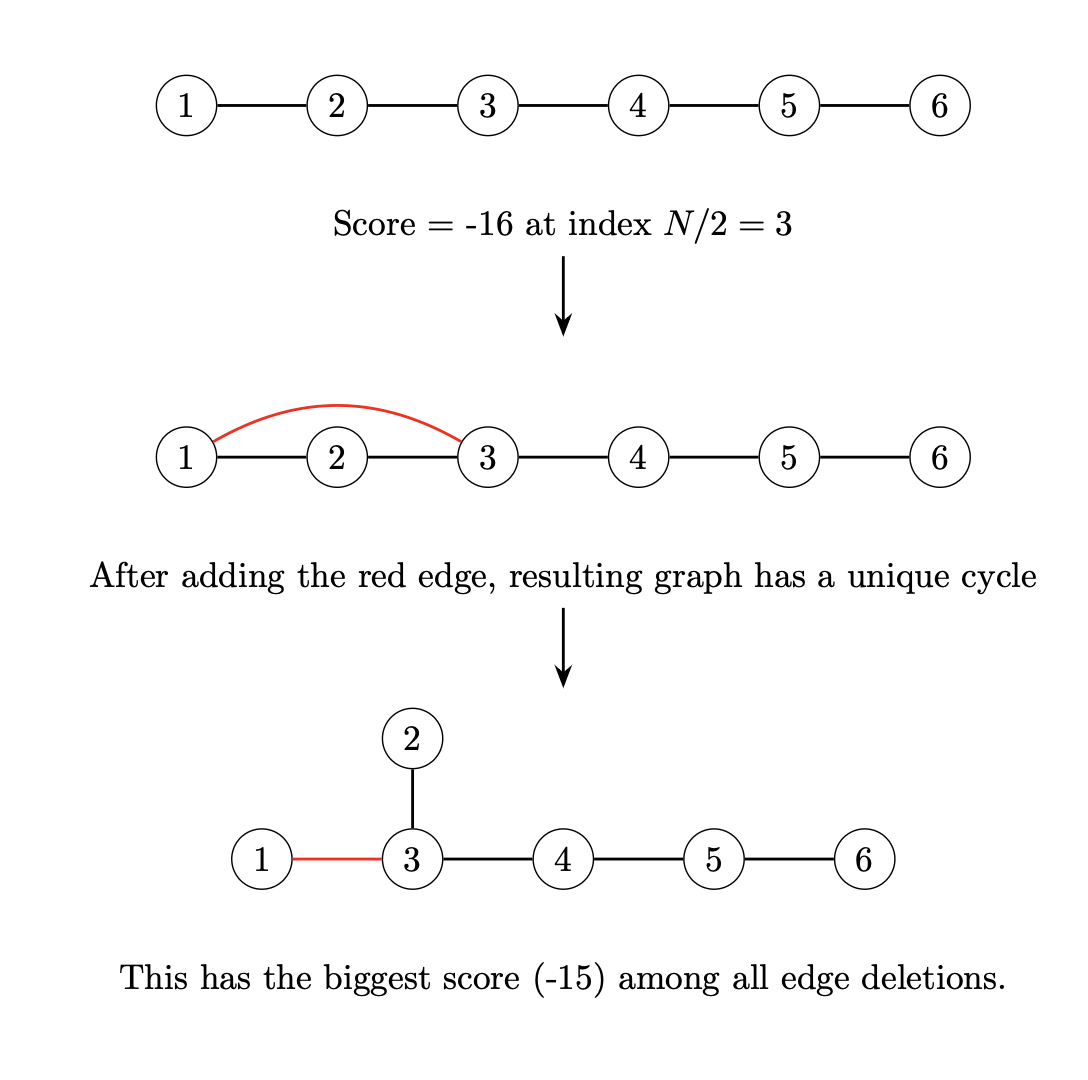}
    \caption{An example of a single edge swap}
    \label{fig:swap}
\end{figure}

In brief, we used a simple edge swap algorithm for our local search. For a given input tree, we first collect all of its \emph{non-edges}. One non-edge in this collection is selected to be added to the tree, creating a unique cycle. We then move along the edges of this cycle, deleting one at a time and checking each time the score of the resulting tree. Among all of the trees created in this way, we replace our original tree with the tree with the highest score. If the deleted edge is not equal to the one that was originally added, then we say that an \textbf{edge swap} has been performed. We continue in this way until a certain number (determined by the user) of swaps have happened. The set of non-edges is determined by the initial input tree, and never updated during the process. This prevents the machine from wasting it time retreading old trees that had been disqualified already. See Figure \ref{fig:swap} for an illustrated example of a single edge swap beginning from the path.

Generally speaking, we found that PatternBoost performed best when the total number of swaps was set to 10. The reasons for this are manifold:

\begin{itemize}
    \item Only performing a single swap was ultimately not a significant enough improvement on the score function for the machine to learn efficiently once the number of vertices grew beyond about 30. In our very initial experiments, we were able to recover the tree $\mathbf{3^*_{34}}$ (see Figure \ref{fig: first counters}) using only a single edge swap, running PatternBoost on our local machines (i.e. without meaningful computing power). Interestingly, it wasn't until we allowed more swaps that we began finding $\mathbf{3_{44}}$.
    \item On the other extreme, performing too many swaps became computationally unwieldy fairly quickly, even on significantly more powerful GPUs. For instance, even considering the smallest interesting case of $N = 26$ vertices, our code running 10 swaps took about an hour and a half to finish 5 epochs. With all of the same parameters, but with each optimization performing the maximum number of edge swaps, the code took almost 5 hours to complete 5 epochs.
    \item Moreover, optimizing too aggressively is -- perhaps unintuitively -- at odds with what we are trying to accomplish here. The transformer at the heart of PatternBoost is designed to learn what high score inputs look like, and from this training set produce representative samples. However, in this project we aren't really interested in finding trees with the highest possible score. We want to produce any and all trees with a positive score. In one of our earlier trials, where the number of swaps was set to its maximum value, we indeed found that the final outputs of the transformer were repeated instances of only two distinct trees. While these were both indeed counter-examples to log-concavity, there were many trees found in the earlier epochs which were lost.
\end{itemize}

\begin{remark}
    As an extra layer of protection against counter-examples being optimized out of existence, our local search algorithm always checks the score of any input tree, and ends if that score is positive. Very importantly, note that this check \emph{only} happens at the initiation of the local search. In other words, if an input tree begins with a negative score, but after some number of swaps becomes positive, this will not cause the search to immediately end.
\end{remark}

Ultimately, after much trial and error we found that 10 swaps was the ``sweet spot" that avoided the troubles of both extremes. Note that, this was only the case when we began running experiments on GPUs. While we were able to run multi-swap experiments on our local machines, this was only for the fewest numbers of vertices, and these experiments often times took many hours to run. It should be mentioned, however, that in local machine experiments on $N = 26$ vertices, we were only able to find $\mathbf{3_{44}}$ when the number of swaps was set much higher. See the discussion before Remark \ref{rmk:hardToFind} for more on this.

It should also be noticeable that we were intentionally vague in our description of how the local search decides which non-edges to add when attempting to perform each edge swap. That is because there are actually many ways to do this, each of which can produce genuinely different behaviors in the machine!

In our code, there is a parameter \orderedges that the user can set to 1, 0, or -1.
\begin{itemize}
    \item If \orderedges is set to 1, the the non-edges are added in decreasing order with respect to the total degrees of their endpoints.
    \item If \orderedges is set to 0, the non-edges are added in lexicographic order. Note that this makes sense, as our vertices are labeled.
    \item if \orderedges is set to -1, the non-edges are added in a random order.
\end{itemize}

\begin{remark}\label{rmk:Charton}
    Upon the recommendation of Francois Charton, we also ran a number of experiments where the total number of swaps was set to 0. This version of the code skips local optimization entirely, and has the transformer try to train itself based exclusively on score feedback. None of these experiments were successful in producing counter-examples.  In all cases, the transformer converged on outputting different permutations of (some minor alteration of) the path graph just after the 10th epoch (see Section \ref{sec:largeEpoch} for more on this particular behavior).
\end{remark}

We think of the three choices presented above as representing a \textbf{sorted} approach to edge addition, an \textbf{unsorted} approach to edge addition, and a \textbf{randomized} approach to edge addition, respectively. While piloting our early experiments on our local machines, we were initially unable to find either of the counter-examples in the $N = 26$ case by using the second or third approaches to edge addition. One thing that became very noticeable, however, is that the machine always ultimately began producing almost exclusively different labelings of the path graph (Once again, see section \ref{sec:largeEpoch}). As such, we created the \orderedges $= 1$ option with the idea that this would promote branching in the trees that were produced during the local search, as non-edges are always added in order to maximize the degrees of their endpoints. Having made this change, we finally were able to find the small counter-examples to log-concavity.

Interestingly, when we got our software onto more powerful GPUs, we discovered that making this choice usually \emph{hindered} our results when $N$ was sufficiently large. For instance, we ran three experiments on $N = 60$ vertices, each of which making a different choice in terms of how edges were added but were otherwise identical in their parameters. After 5 epochs:

\begin{itemize}
    \item the \orderedges $= 1$ (sorted) program had found 247 counter-examples to log-concavity;
    \item the \orderedges $= 0$ (unsorted) program had found 26,766 counter-examples to log-concavity;
    \item the \orderedges $= -1$ (randomized) program had found 24,163.
\end{itemize}

While the performance of the randomized and unsorted programs were roughly equivalent, the sorted program was astonishingly worse. Our ultimate takeaway from this is that with the computational limitations of working on our personal computers, the transformer was not provided enough data to produce good results. As such, we had to modify the local search in such a way to encourage the machine in directions that were more likely to yield the trees we were hoping for. As we pushed things to scale, however, doing this was ultimately more harmful than good. It seems that counter-examples to log-concavity display an incredible diversity of different shapes and behaviors, and so the best way to encourage the discovery of these trees was to not bias the machine one way or the other.

All of the above being said, however, in our experiments at $N = 37$, only the ordered version of the local search was able to find any trees with breakage at $N/2$.

Finally, we have already noted that our experiments, if run long enough without finding any positive scores, would frequently converge on producing mostly permutations of (slight alterations of) the path graph. While we will go into more details about this behavior in Section \ref{sec:largeEpoch}, it is worthwhile to note here that in response to this we created a version of the local search which explicitly punished the path. Specifically, if the \prufer code of a path graph was entered into the local search, the search first replaced it with a star graph. The star graph was selected for a number of reasons. Firstly, from a coding perspective, it is simple to replace the origin \prufer code with the code $[1,1,\ldots,1]$, corresponding to the star with center 1. But more importantly, the independence number of the star graph is $N-1$. It turns out, mostly because $N/2$ lies in the middle of the independence sequence rather than near the ends, star graphs are trees whose scores at indices near $N/2$ are colossally negative relative to other trees.

Through experimentation we have discovered that this "path punishing" version of the local search did not perform meaningfully better when looking for breakage at $N/2$. However, as we will see in section \ref{sec:N/2-1}, this change was the final key necessary to begin producing examples of breakage at index $N/2-1$ when $N$ was even.

\section{Notable results and notable failures}

\subsection{Notes on parameters}\label{sec:parameters}

Before we dive into the various results and non-results, we take a moment to discuss the various parameters involved in running our experiments. In Section \ref{sec:local}, we mentioned the importance of choosing the right number of swaps to perform during the local optimization phase, as well as the different behaviors that manifest when choosing how to sort the added edges. Rather than discuss those local parameters again, we instead take this time to discuss the global parameters related with the behavior of the transformer and its training.

To begin, we spent a lot of our time running experiments on our local machines. During these experiments we generally:

\begin{itemize}
    \item Seeded the process with 50,000 \prufer codes generated uniformly at random.
    \item At the end of each local phase, the top 50,000 scoring trees were send to the transformer. Of these, 49,000 were used for training, while the remaining 1,000 were used for testing. 
    \item Training occurred in a loop over 1,000 iterations. During these iterations, the training and testing loss never diverged. Both gradually decreased over the course of the entire PatternBoost trial.
    \item After completing training, we queried the transformer for 2,000 samples. Those samples corresponding to legitimate \prufer codes were combined with the original 50,000 training and testing codes, and sent back to the local phase.
    \item This continued for around 3 epochs.
\end{itemize}

Obviously, the numbers given above are nowhere close to scale, and indeed our experiments were ultimately only able to find counter-examples in the cases of $N = 26, 28,$ and $30$ vertices. Beyond 30 vertices, the computation became too heavy for our machines. Further note that we were unable to find counter-examples in the odd cases $N = 27$ and $29$ during these early attempts. In fact, as of the writing of this article, we have still never been able to use our methods to produce counter-examples with breakage at $N/2-1$ when $N$ is odd. This suggests that the odd case is fundamentally more difficult for the machine to wrap its head around. We will provide what we believe is the reason for this in Section \ref{sec:alpha}.

Following these early excursions, we were finally able to get access to GPU cores thanks to the combined efforts of the Stevens Institute for Artificial Intelligence and The Stevens Lab for AI in Mathematics Education. 

Our experiments were run on Nvidia RTX A6000 GPUs. We generally ran 64 concurrent threads for each experiment. During these experiments we usually:

\begin{itemize}
    \item Seeded the process with 50,000 \prufer codes generated uniformly at random, as on our personal machines.
    \item At the end of each local search, the top 50,000 scoring trees were send to the transformer. Of these, 49,000 were used for training, while the remaining 1,000 were used for testing. Once again this wasn't changed from the previous setting. 
    \item Training occurred in a loop over 8,000 iterations. As before, testing and training loss maintained a close distance, usually only differing after the second decimal point.
    \item After completing training, we queried the transformer for 100,000 samples. Those samples corresponding to legitimate \prufer codes were combined with the original 50,000 training and testing codes, and sent back to the local phase.
    \item This continued between 5 and 15 epochs.
\end{itemize}

A few remarks are necessary with regards to the parameters discussed above. Firstly, the choice of 8,000 training rounds was made after experimentation. Our very first GPU experiments used 12,000 training rounds, but it was noted that there were diminishing returns. Generally, testing and training loss decreased meaningfully until around 8,000 training rounds. Beyond this point improvements in testing and training loss were much less common.

Insofar as our sampling step, it is notable that the machine was generally quite good at not producing non-sense. For instance, in practice, in each epoch some number of \prufer codes at or above 148,000 were preserved and sent to the local phase. We believe that the reason for this resiliency to hallucination is precisely because of how nice \prufer codes are for transformers. Tokenization is one to one in this case. Every token can be inserted into the string at any point, and the length of the tokenized \prufer code is always the same as the length of the original \prufer code.

Finally, one may look at the number of epochs and feel that it is quite lower than one expects. In cases where we looked for counter-examples whose log-concavity breakage happened at the index $N/2$, the machine generally started producing counter-examples after the third generation at the latest. In fact, for values of $N$ near the smallest non-trivial value of 26, counter-examples were found within the samples produced by the transformer after only one epoch. In all cases, after around the 5th epoch, there were diminishing returns in terms of the number of counter-examples that were being produced. For instance, we ran one experiment with the above parameters at $N = 60$ for 10 epochs. After the 5th epoch the machine had found 24,163 different counter-examples with breakage at $N/2$. By the end of the experiment, a total of 38,367 counter-examples had been found. In other words, after 5 epochs over 60\% of the counter-examples had been found. All this being said, after creating the version of the local search that punished the machine for producing the path graph, all successful experiments looking for counter-examples whose breakage happens at $N/2-1$ took around 10 epochs before finding their first.

In the next section, we will focus primarily on the things that we observed when we tried to run experiments that lasted for greater numbers of epochs.

\subsection{Behaviors in the large epoch limit}\label{sec:largeEpoch}

As stated in the prior section, our experience shows that simply increasing the number of model parameters does not guarantee better performance. In the $N=60$ case, we saw that the fifth generation marked the point of greatest production, after which counter-examples were harder to come by. This obviously suggests that the learning hit some bottleneck or "stable state". There are a number of possible explanations for this, including:

\begin{itemize}
\item The security measure, discussed in Section \ref{sec:local}, which skips optimization for trees whose score is already positive. It is possible that there are counter-examples to log-concavity which are most easily obtained through edge-swaps performed on other counter-examples. This might have had a small biasing effect on the machine;
\item After around 5 epochs, the machine has found the majority of all possible counter-examples to log-concavity with breakage at whichever index is being targeted. While it is possible that this is the case, we find it unlikely considering the behaviors that we observed in the independence numbers of the samples being produced by the machine. We will go into more detail about this in Section \ref{sec:alpha};
\item Some other aspect of our chosen local search is biasing the machine into only finding a particular subset of all possible counter-examples. While this is definitely possible, we do not have evidence one way or another to say. The remarks made in Section \ref{sec:local} on the difference in performance between ordering the non-edges or not should make it clear that this is a genuine concern.
\end{itemize}

One behavior that is very interesting to note here is that the machine seemed to have a hard time rediscovering counter-examples. What we mean by this is that when you look at the final distributions of any experiment that successfully found counter-examples, positive scores would by and large each only appear once, and occasionally twice. One thing that this implies is that the machine never quite figured out that it could repeat positive scores by permuting the labels in the tree. Note that, looking over the counter-examples in the Appendix of this work, counter-examples do not seem to be particularly symmetric and so you would expect there to be many different \prufer codes which correspond to the same tree. Importantly, this is in stark contrast to trees like the path, which the machine would frequently become intensely focused on and produce thousands of repeats.

The above being said, in experiments where we looked at small vertex sizes, such as $N = 26$, and allowed the machine to run for a sufficiently long time, it did eventually start repeating counter-examples more freely. In one experiment with $N = 26$, and performing the maximum amount of edge swaps (rather than just 10), after 6 epochs the 50,000 best scoring outputs of the machine included 4390 copies of $\mathbf{3_{44}}$, 96 copies of $\mathbf{3^*_{34}}$, and 45514 copies of the path.

\begin{remark}\label{rmk:hardToFind}
    As an interesting side note, it was a repeated theme in our work that $\mathbf{3^*_{34}}$ was ``easier to find" than $\mathbf{3_{44}}$, in the sense that it took us many more experiments on our personal computers before we were able to find $\mathbf{3_{44}}$. However, as the above results of our exhaustive super computer search show, once the machine found $\mathbf{3_{44}}$, it had a far easier time finding its various label permutations. We are unsure why this is. It is possible it is related with the fact that $\mathbf{3_{44}}$ has a strictly higher score than $\mathbf{3^*_{34}}$,
    
    In fact, it is noticeably unintuitive that this would be the case because $\mathbf{3^*_{34}}$ is clearly less symmetric than $\mathbf{3_{44}}$, which implies there are more \prufer codes that correspond to it.
\end{remark}

All of what has been reported thus far has related with the long term behaviors of the machine during successful experiments. Throughout our explorations, once we were able to run things on GPUs, we were never unsuccessful at finding counter-examples with breakage at index $N/2$. These experiments spanned a large number of the vertex sizes from 26 to 101. However, for a very long time our methods were entirely unable to find counter-examples with breakage at any index other than $N/2$. In fact, as of writing, we have only partially overcome this shortcoming by using the machine to find counter-examples whose breakage happens at $N/2-1$. We will discuss that case in more detail in Section \ref{sec:N/2-1}.

\begin{figure}
    \centering
    \includegraphics[width=0.35\linewidth]{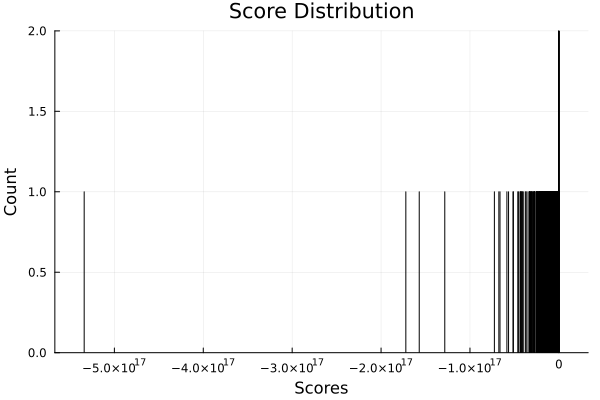}\includegraphics[width=0.35\linewidth]{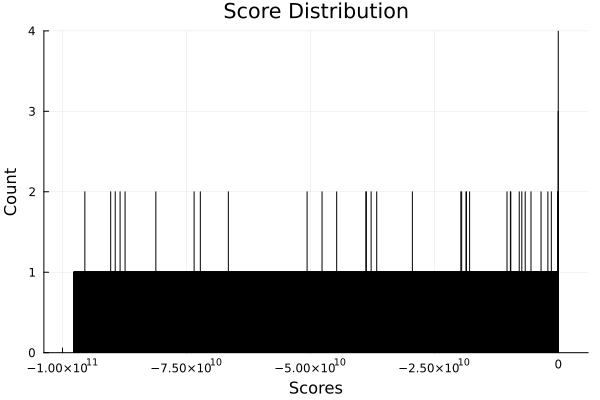}\includegraphics[width=0.35\linewidth]{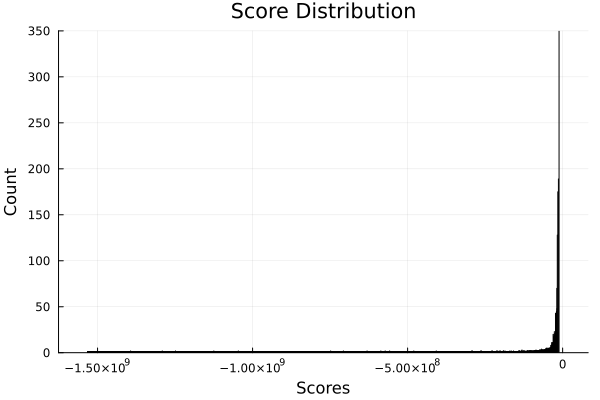}
    \includegraphics[width=0.4\linewidth]{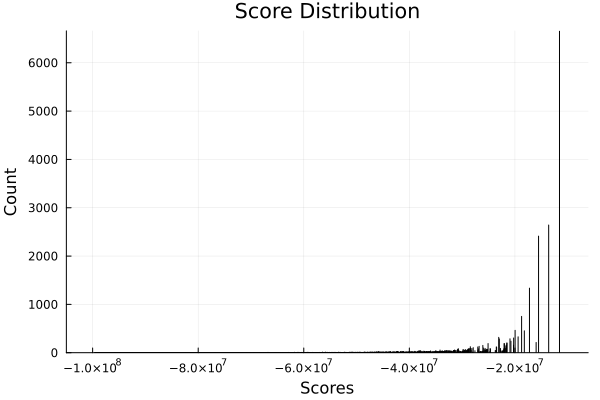}\includegraphics[width=0.4\linewidth]{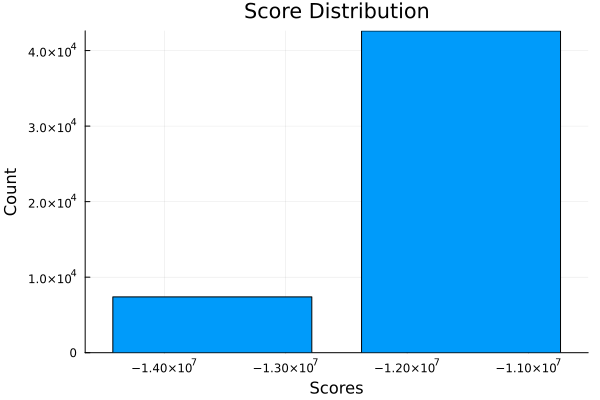}
    \caption{The top 50,000 scores in each epoch during our first attempts to find breakage at $N/2-1$. Here, $N = 56$}
    \label{fig:eventualTree}
\end{figure}

Before diving into that, however, we do want to say some words about the things that the machine was doing in these cases. For instance, from the work of Galvin \cite{galvin2025trees}, it was known that there should be at least one counter-example with $N = 56$ vertices whose breakage happens at index $N/2-1 = 27$. Early on we tried to find this example using our methods. The histogram plots of the top 50,000 scoring trees across all 5 epochs can be found in Figure \ref{fig:eventualTree}.

The tall column near 0, on the right most side of each plot, ultimately becomes one of only two columns by generation 5. This column corresponds to the various relabeling of the path graph. To be more precise, the distribution of scores after the 5th generation is 42608 copies of the path, and 7392 copies of a path graph on 55 vertices, with an extra leaf appended to one of its inner vertices (similar to the final tree in Figure \ref{fig:swap}).

The main takeaway, was that in almost all cases we attempted, if the machine ran for long enough without finding a counter-example, it would somehow always default to producing exclusively the path or path-like trees. It remains somewhat of a mystery to us why it does this. For instance, the path is \emph{not} a local maximum of our local search algorithm until $N$ is sufficiently large relative to the index being checked. For instance, when the score is relative to index $N/2$, the path does not being a local extrema until $N = 60$. Despite this fact, this ``convergence to the path" behavior is consistent in all experiments we ran, both at values of $N$ before and after the threshold for when the path becomes a local optimum. Our best guess for why this happens is that the path is the simplest example of a tree with the smallest possible independence number, $\alpha = N/2$ when $N$ is even, and $\alpha = N/2+1$ when $N$ is odd. Moreover, a \prufer code corresponds to a path if and only if no entry repeats. This ease of creation might also influence the machine in some way.

We will speak more on the independence number in the next section.

Once we implemented the version of the local search that explicitly punishes the path, this tendency to converge on a single tree largely stopped.

\subsection{The curious behavior of the independence number}\label{sec:alpha}

One of the most interesting behaviors in the samples produced by the transformer is that they tended to have extremely small independence number. It is a classical result that the independence number of a uniformly randomly generated \prufer code of length $N-2$, has, in the large $N$ limit, an expected value around $.57N$ \cite{meir1973expected}. Despite this fact, for all $N$ up to at least 101, when looking for breakage at $N/2$ the graphs produced by the transformer almost all have independence number equal to $N/2 + 1$. In fact, when $N$ was even, the machine almost never produced samples whose independence number was strictly above $N/2+1$. It was only in the cases where $N$ is odd that the machine sometimes would produce samples with independence number equal to $N/2+2$, but never higher.

This behavior is interesting for at least two reasons:
\begin{itemize}
    \item We do not see anything in our local optimization step that explicitly encourages the independence number to be as small as possible. Our best guess for why this is happening is because of our score function. Elements in the independence sequence will always tend to be much smaller near the top of the sequence, and so if we are looking at the coefficients around $N/2$, it is likely helpful to force these numbers to be as small as possible by shrinking the independence number.
    \item As far as the authors are aware, there is no obvious way to read the independence number off the \prufer code. It is therefore somewhat shocking how proficient the machine gets at doing this. In most of our experiments looking for breakage at $N/2$, the independence number is forced down to $N/2+1$ in no more than one or two generations.
\end{itemize}

Recall in previous sections we discussed the tendency of the machine to converge on the path if given enough time without finding positive scores. It is our belief that this tendency for the machine to shrink the independence number is a big part of what is causing this. Rather than the path being a local extrema with respect to our optimization algorithm, it is an absolute extrema with respect to the independence number. This is also our best theory for why the machine had so much of a harder time finding breakage when $N$ was odd. Because we decided to define $N/2$ in the odd case to be the floor rather than the ceiling, where we are searching for breakage will always be further away from the minimum independence number than it is in the even case.

This behavior also explains why Galvin's constructions rarely appear in our database of counter-examples \cite{galvin2025trees} (Though we did find the tree $T_{4,4}$, as depicted in Example \ref{ex:T44} of the Appendix). In Galvin's construction of the family $\{T_{m,t}\}$, each tree in the family has a total of $1+m+2mt$ vertices, with $\alpha(T_{m,t}) = (1+t)m$ the log-concavity of the independence sequence breaks at $mt+2$. In particular, the independence number can be made arbitrarily bigger than the minimum of $N/2$, and the gap from the independence number to the index of breakage also grows without bound. The current machine, however, focuses on counter-examples with small independence numbers.

When we ran our (successful) trial on $N=56$ vertices, which has the same number of vertices as $T_{5,5}$ in Galvin's constructions, we found 6 trees whose independence sequence was not log-concave, but none of them was $T_{5,5}$.

\begin{conjecture}
    If $T$ is a tree on $N$ vertices whose independence sequence is not log concave in index $N/2 - \beta$, then the independence number of $T$ must be at least $N/2+\beta + 1$.
\end{conjecture}

We came to this conjecture through our various experiments. Noting all of the things discussed above, if this conjecture were true, it implies that our machine will have a harder and harder time finding counter-examples as we look further down the sequence. 

\begin{remark}\label{rmk:smallBanFail}
Considering our success of blocking the path graph to find breakage at $N/2-1$, one might expect, along with the above conjecture, that it may make sense to write a local search that punishes all trees whose independence number isn't high enough. While we have coded such a local search, it has been disappointing in terms of its results. Rather than forcing the machine to focus on building trees with higher independence numbers, in practice it seems to just cause the machine to become confused. In the few experiments we have run with this code, scores do not seem to meaningfully improve. It is possible that the ability to shrink the independence number is so fundamentally important to the way the machine learns to improve scores, that preventing it from doing so causes the machine to perform considerably worse.
\end{remark}

\subsection{Searching away from half the number of vertices}\label{sec:N/2-1}

In our efforts to search away from half the number of vertices, the vast majority of our trials have ultimately ended in failure. While this may be because trees with breakage away from $N/2$ are exceedingly rare compared to those with breakage at $N/2$ on the nose, because of Galvin's work \cite{galvin2025trees} we knew such trees must exist. As such we began focusing our attention on looking for breakage at index $N/2 -1 = 27$ for trees with $56$ vertices, as the smallest tree in Galvin's work with breakage at $N/2-1$ has this many vertices.

We discussed in Section \ref{sec:largeEpoch} that our initial attempts to find these counter-examples failed because the machine repeatedly insisted upon mostly outputting the path. This led us to develop the version of local search that intentionally punishes the machine if it produces the path. Using this new local search we were able to find counter-examples with breakage at index $N/2-1$ when $N=56$, and $58$. Note that we have still have not any success in finding breakage at $N/2-1$ whenever $N$ is odd.

\begin{figure}
    \centering
    \includegraphics[width=0.5\linewidth]{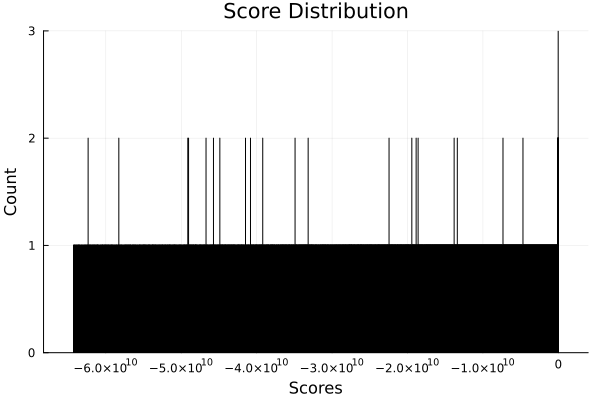}\includegraphics[width=0.5\linewidth]{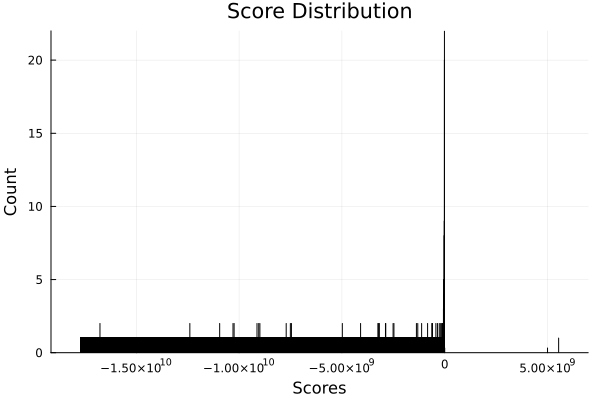}
    
    \caption{The histogram of top scores in generation 6, and generation 7, respectively, when looking for log-concavity breakage at index 27 for trees of order 56. The histogram for the 6th generation is quite similar looking to those for the 2nd through 5th generations as well.}
    \label{fig:suddenPositive}
\end{figure}

One thing that is notable about these experiments is how suddenly they found their counter-examples. When looking for breakage at $N/2$, you could see the machine generation over generation producing better scores until finally finding one that was positive. When looking for breakage at $N/2-1$, however, the machine seemed to stumble around extremely negative values (between about $-10^{11}$ and $-10^7$) until the 7th generation when it suddenly found a single extremely large positive score ($\sim 10^9$). You can see an example of this in the plots of Figure \ref{fig:suddenPositive}

The second notable point is that we feel that these experiments truly illustrate the full power of the PatternBoost architecture. In other words, in order to produce these counter-examples \emph{both} the local and global phases of the process were necessary. As noted in Remark \ref{rmk:Charton}, when we ran experiment without any local search the machine never came close to finding counter-examples, and in fact usually converged on something path-like. Conversely however, our original attempts to find counter-examples with breakage at index 27 when $N = 56$ all failed until we began punishing the path. This illustrates that the local search by itself was never going to be able to find these counter-examples, it needed better (i.e. non-path-like) inputs from the transformer!

\printbibliography

\section{Appendix: A random sampling of trees without a log concave independence polynomial.}
\begin{example} A tree with 101 vertices.
    \begin{figure}[H]
    \centering
    \includegraphics[width=\linewidth]{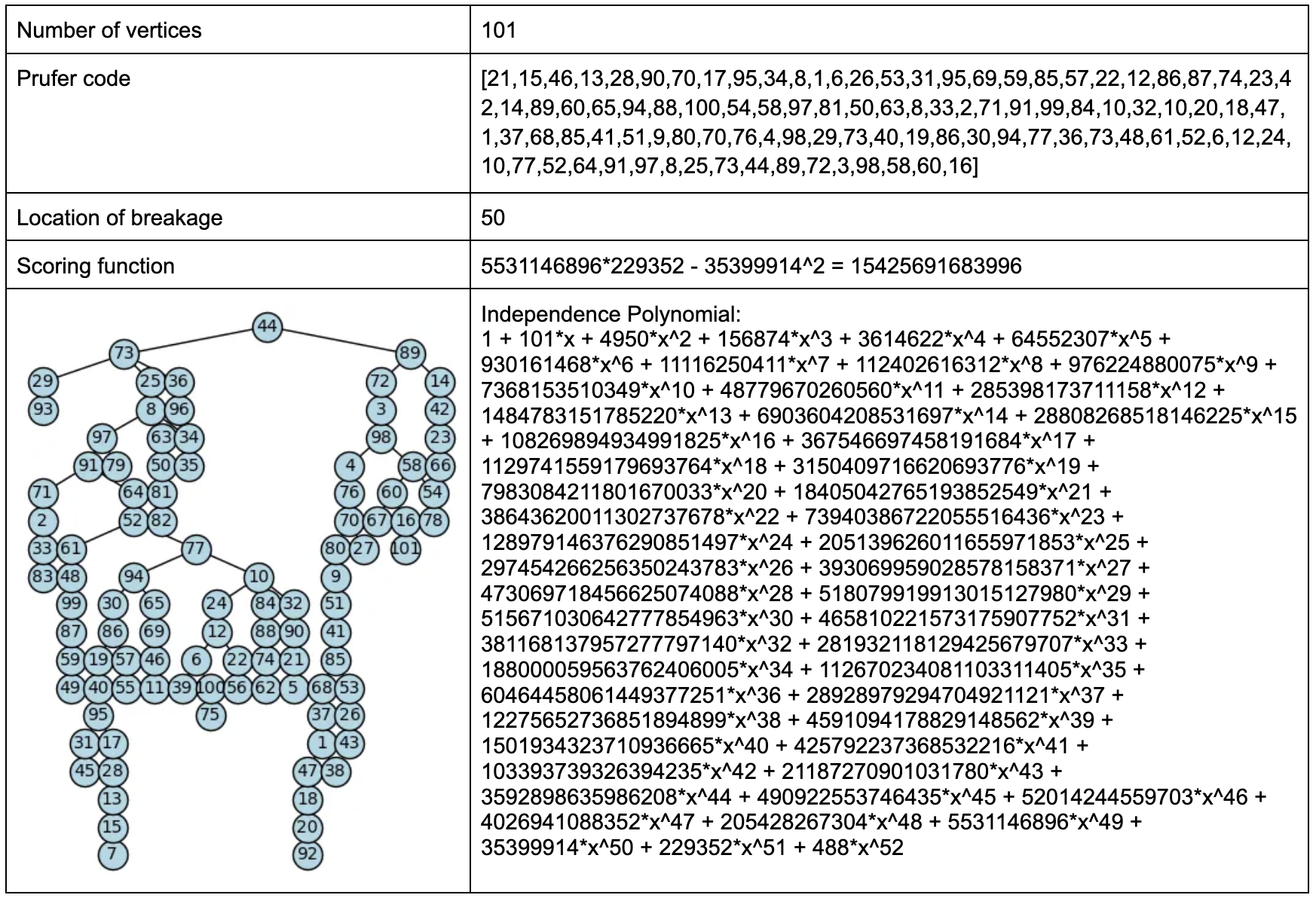}
    \end{figure}
\end{example}

\begin{example} A tree with 67 vertices, with breakage happening at $\alpha -2$.
    \begin{figure}[H]
    \centering
    \includegraphics[width=\linewidth]{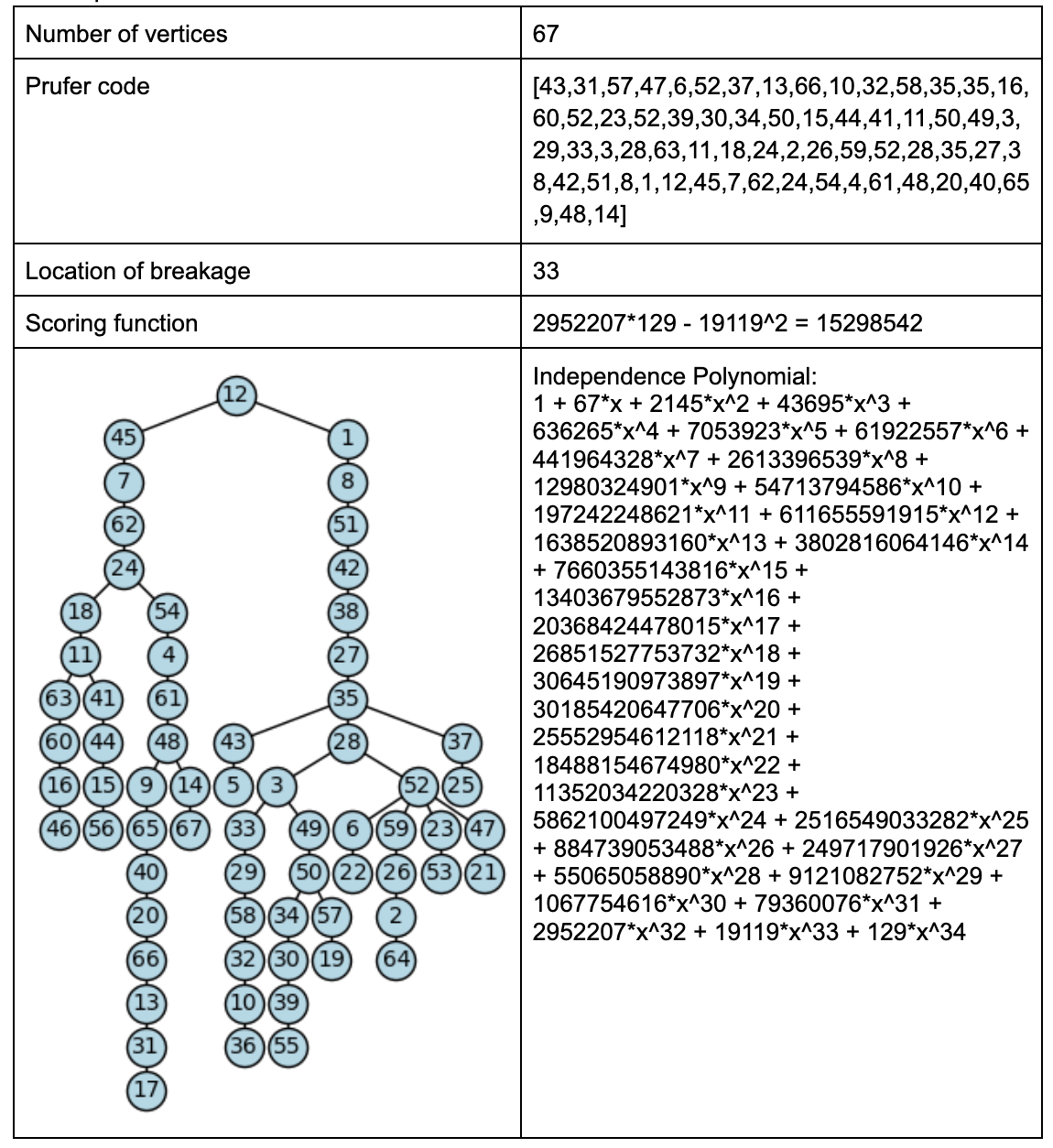}
    \end{figure}
\end{example}

\begin{example} A tree with 67 vertices, with breakage happening at $\alpha -1$.
    \begin{figure}[H]
    \centering
    \includegraphics[width=\linewidth]{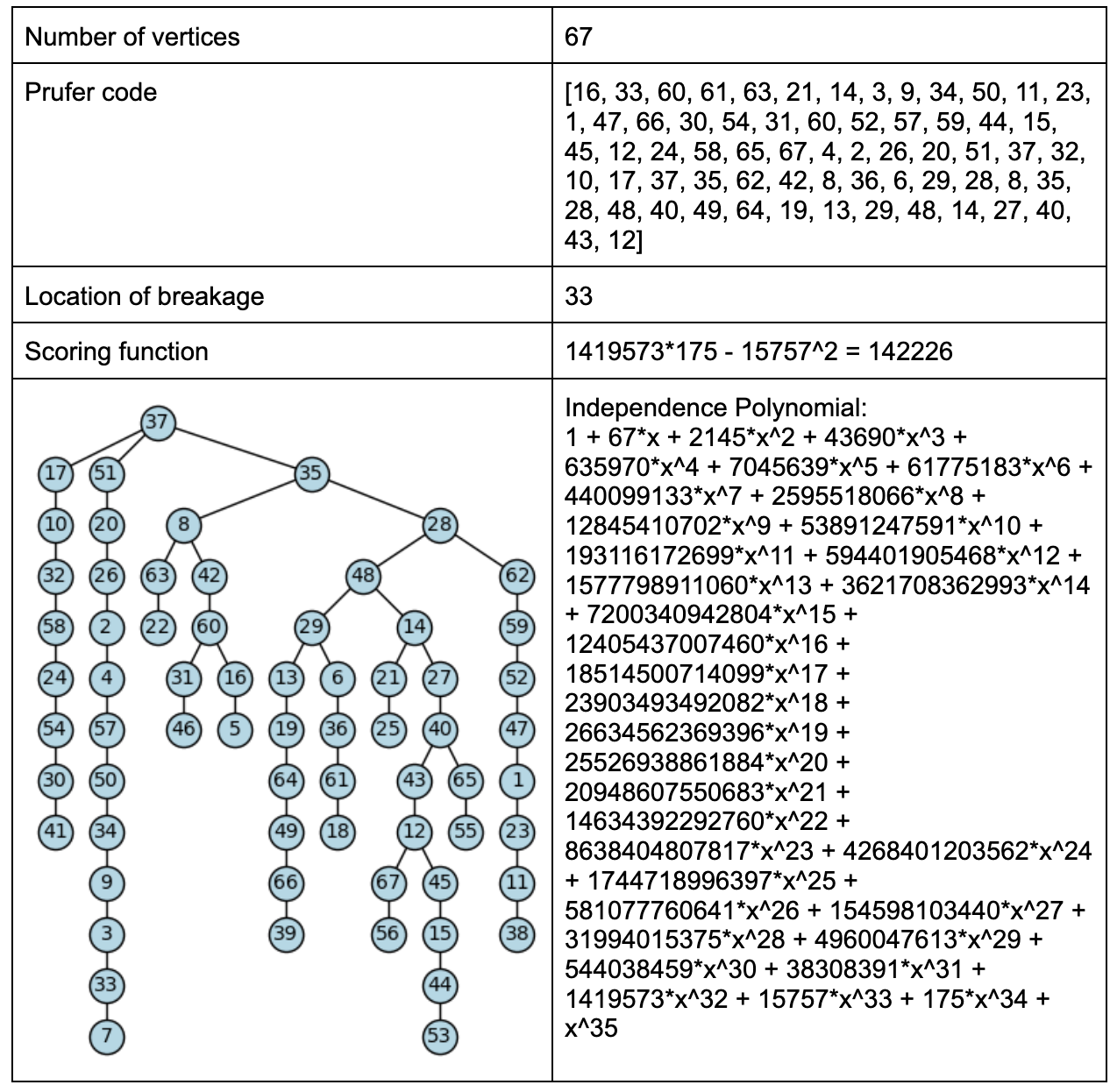}
    \end{figure}
\end{example}

\begin{example} A tree with 60 vertices.
    \begin{figure}[H]
    \centering
    \includegraphics[width=\linewidth]{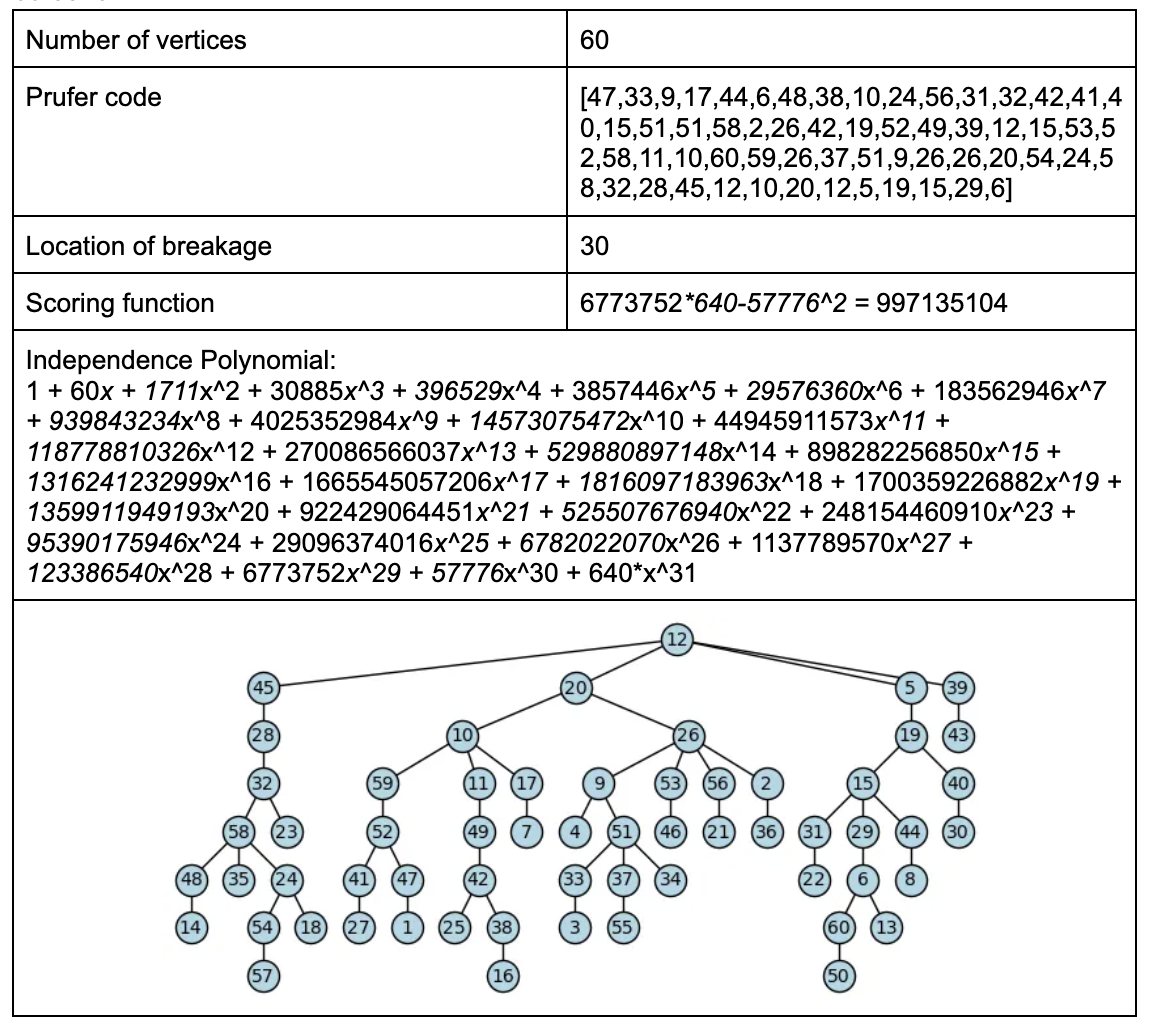}
    \end{figure}
\end{example} 

\begin{example} A tree with 60 vertices, with score = 9.
   \begin{figure}[H]
    \centering
    \includegraphics[width=\linewidth]{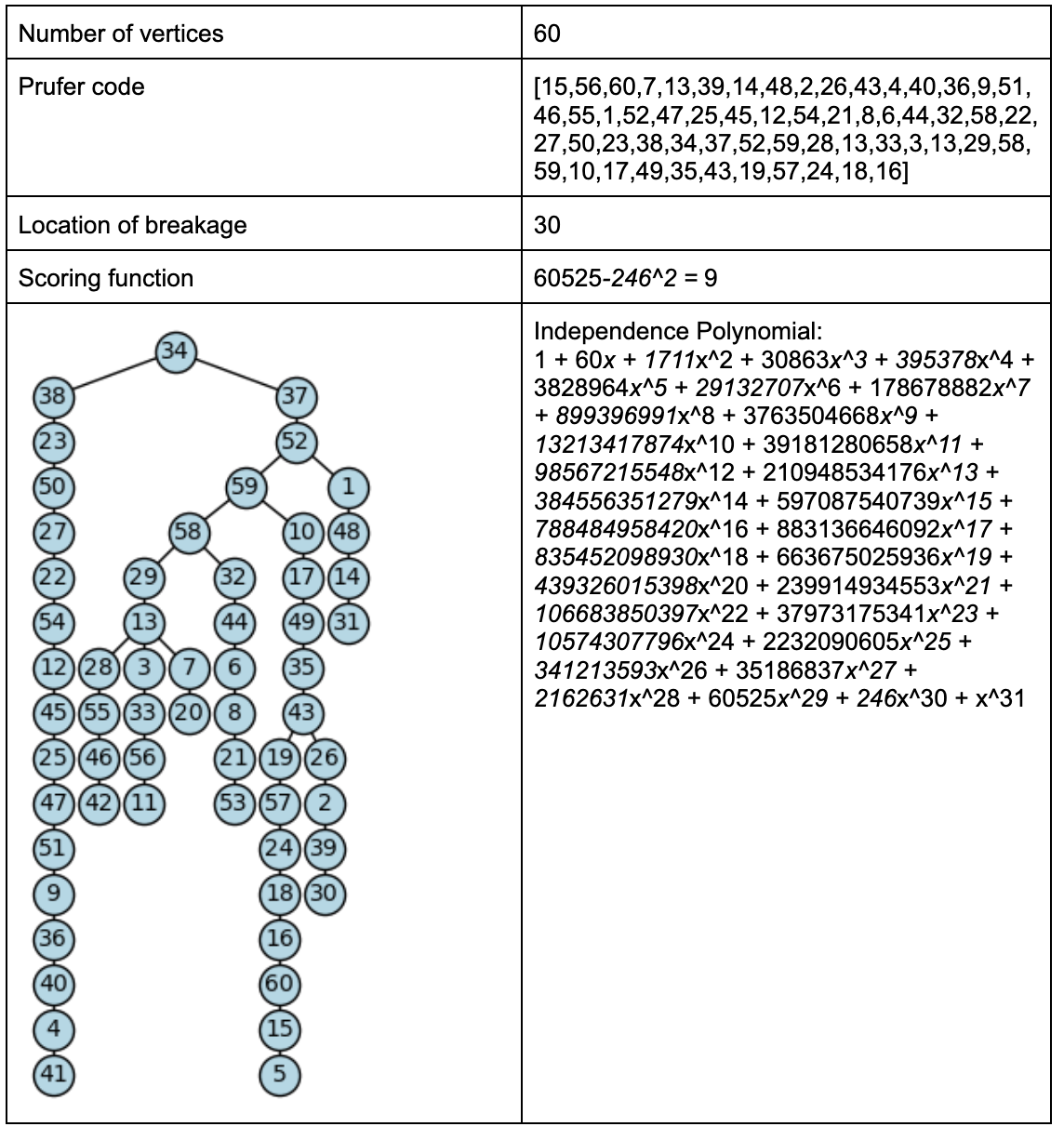}
\end{figure}
\end{example}

\begin{example} A tree with 58 vertices; log-concavity breaks at $N/2-1$.
   \begin{figure}[H]
    \centering
    \includegraphics[width=\linewidth]{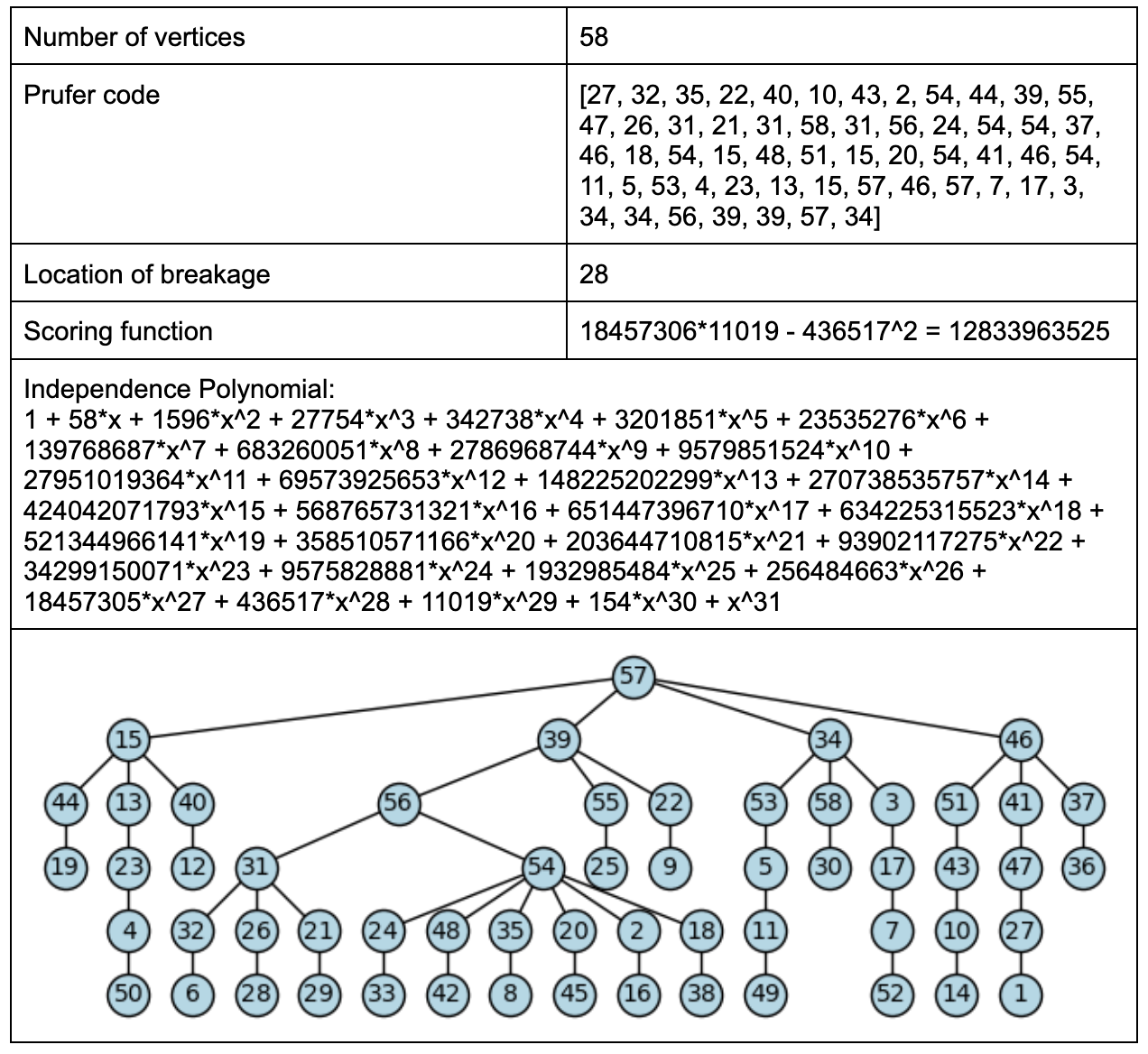}
\end{figure}
\end{example}\label{ex:alpha-3}

\begin{example} Another tree with 58 vertices; log-concavity breaks at $N/2-1$.
   \begin{figure}[H]
    \centering
    \includegraphics[width=\linewidth]{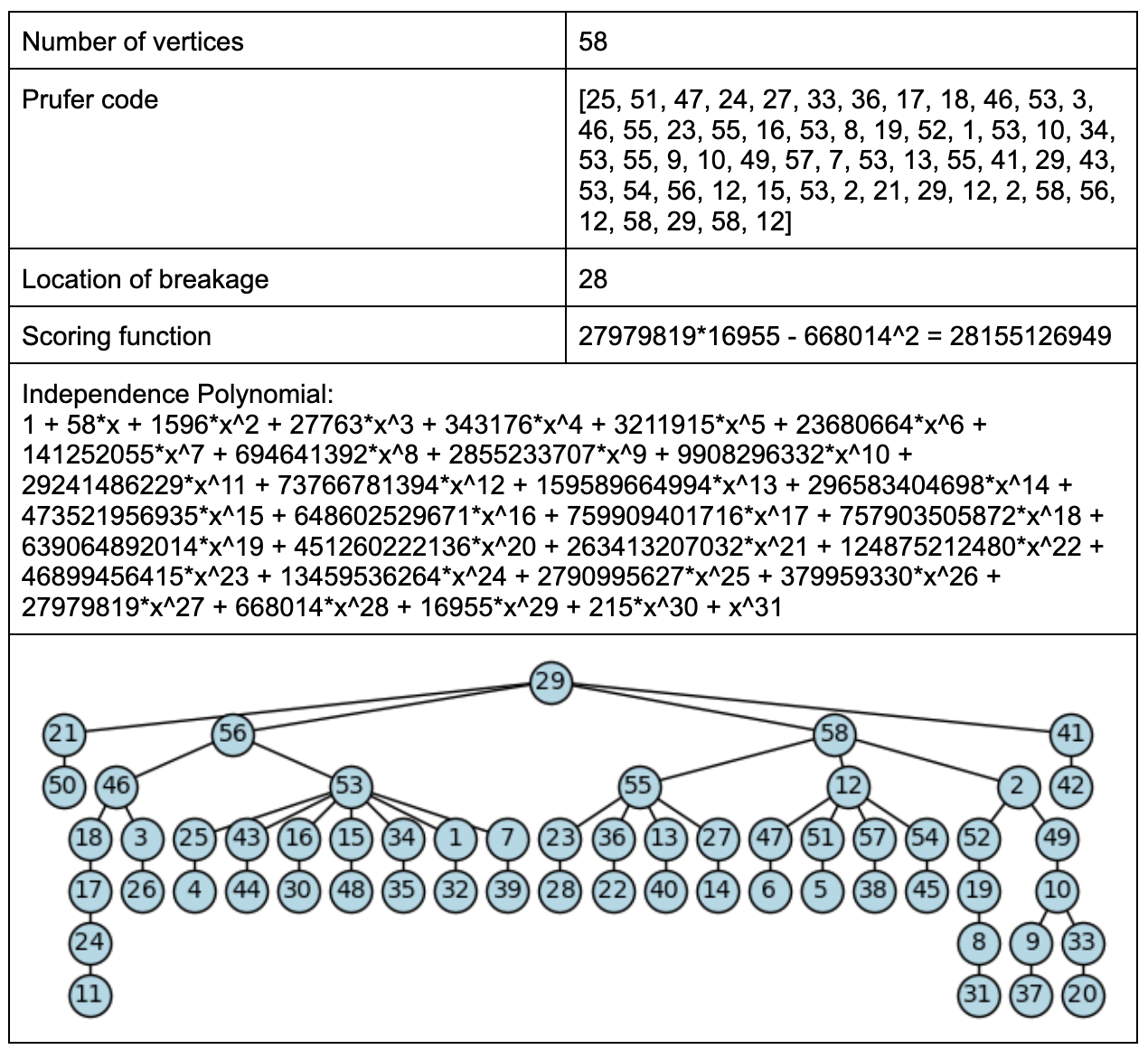}
\end{figure}
\end{example}

\begin{example} A tree with 56 vertices, log-concavity breaks at $N/2-1$.
    \begin{figure}[H]
    \centering
    \includegraphics[width=\linewidth]{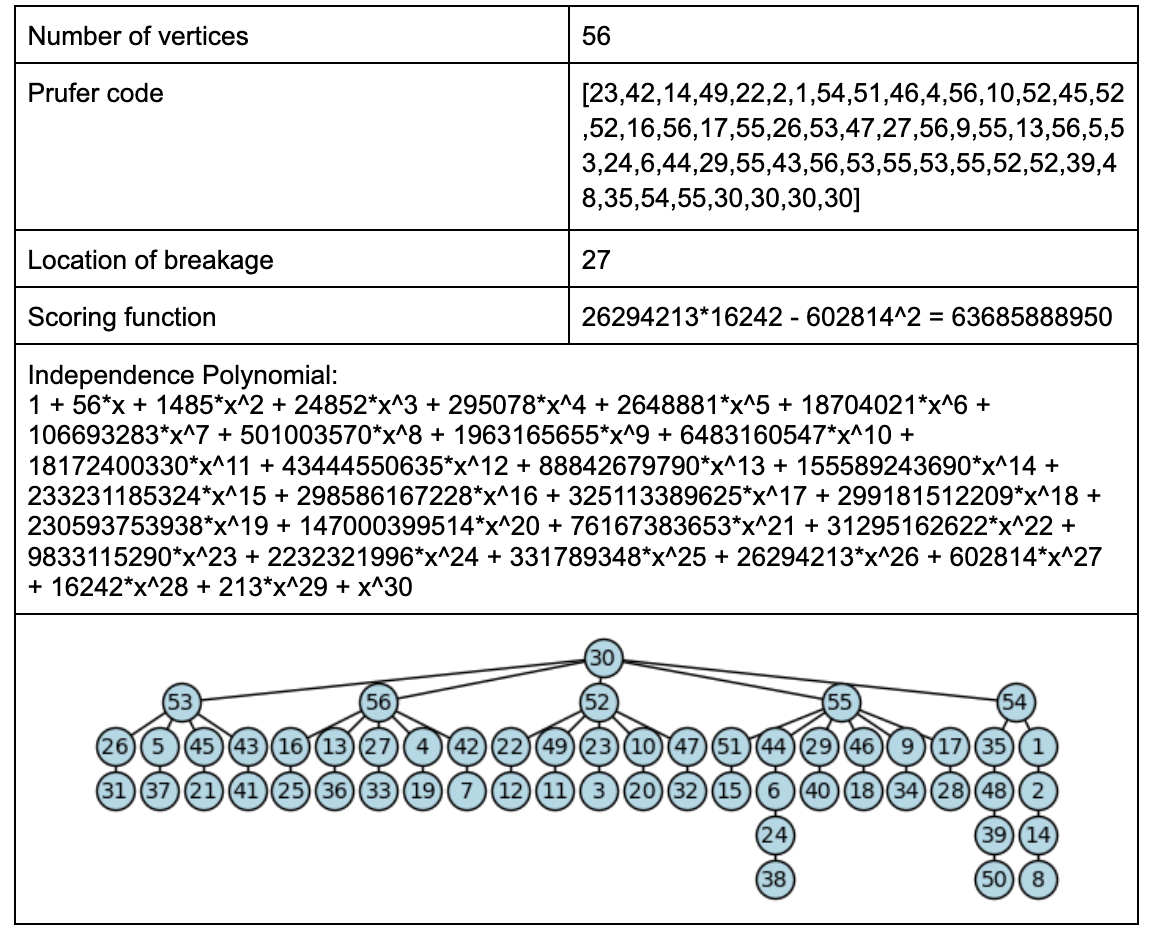}
    \end{figure}
\end{example}

\begin{example} A tree with 56 vertices, with score = 1.
    \begin{figure}[H]
    \centering
    \includegraphics[width=\linewidth]{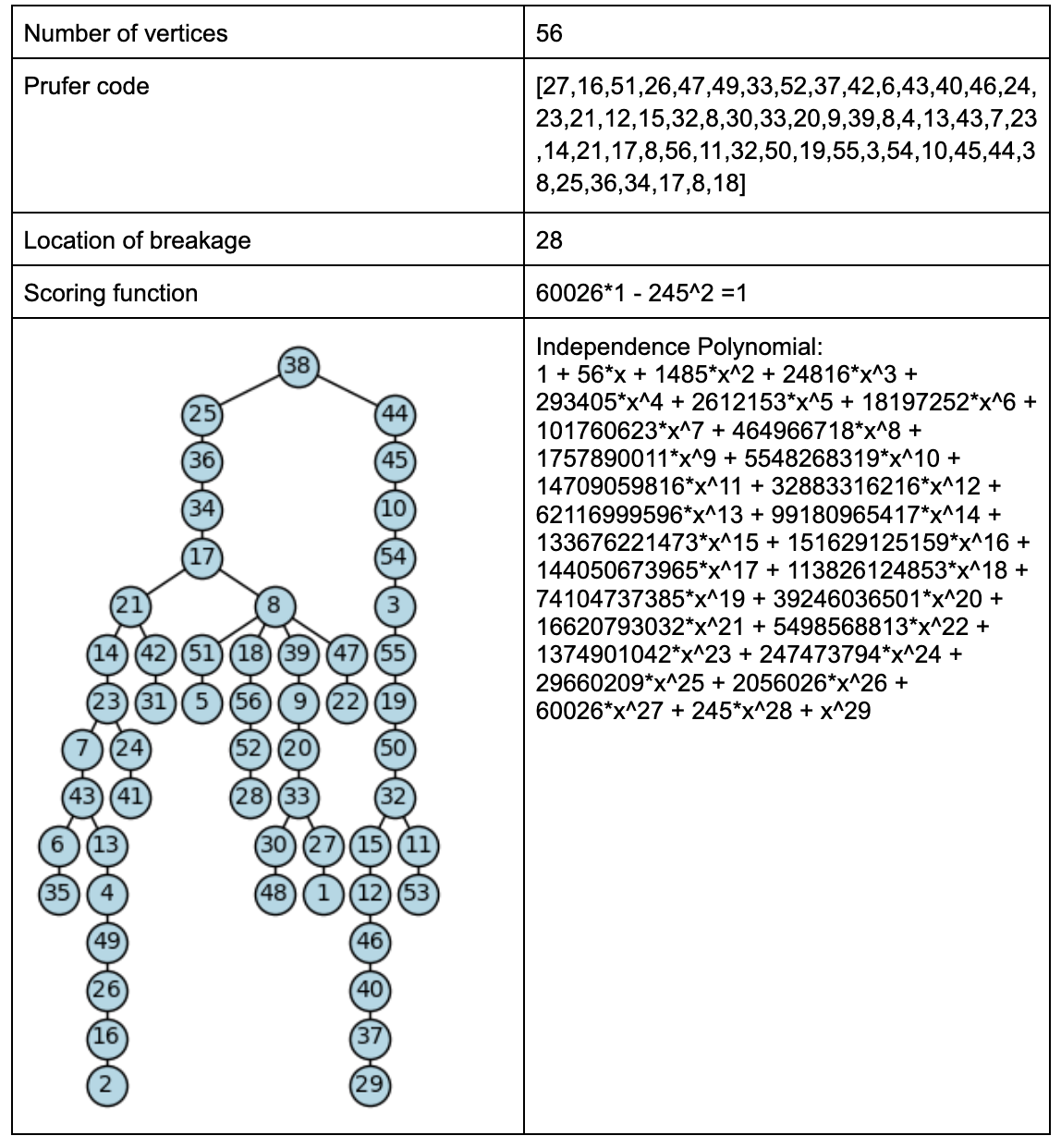}
    \end{figure}
\end{example}

\begin{example} A tree with 54 vertices.
    \begin{figure}[H]
    \centering
    \includegraphics[width=\linewidth]{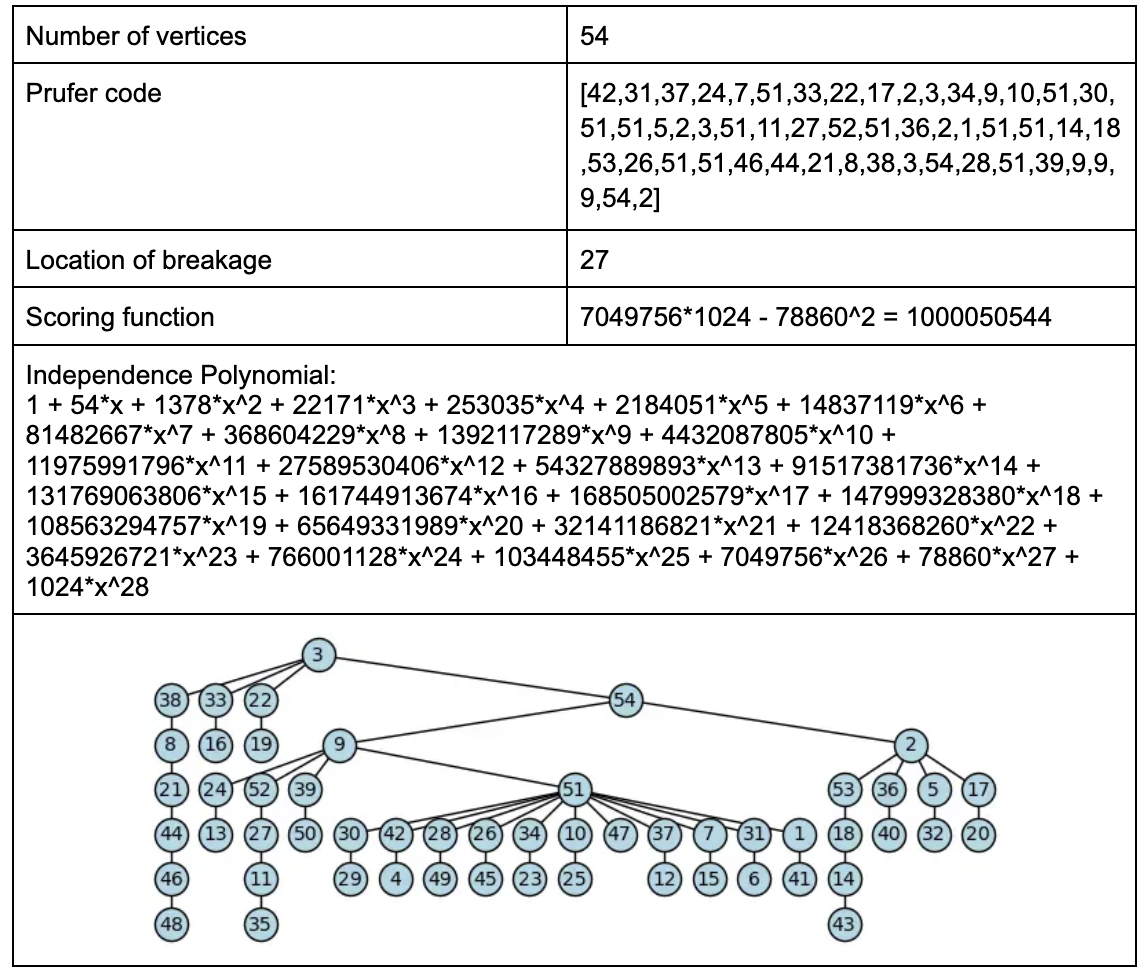}
    \end{figure}
\end{example}

\begin{example} Another tree with 54 vertices.
    \begin{figure}[H]
    \centering
    \includegraphics[width=\linewidth]{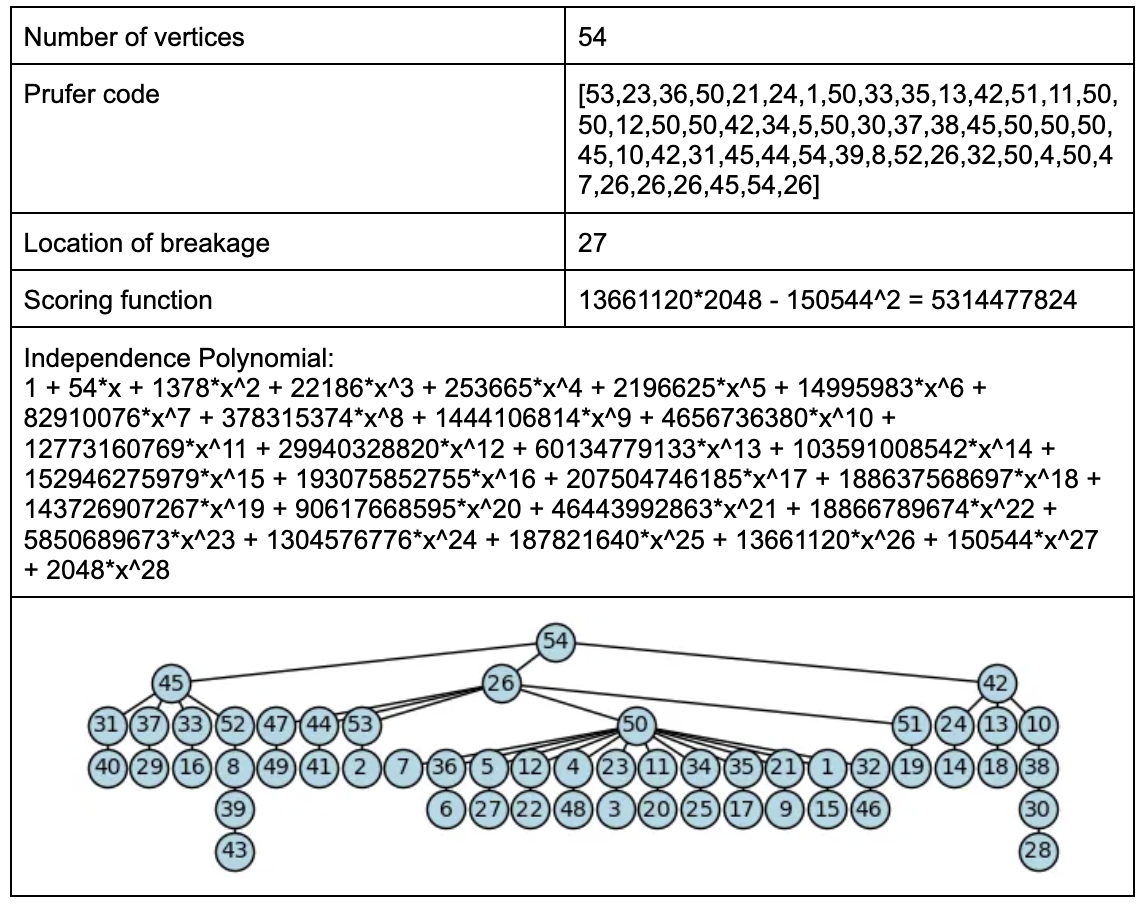}
    \end{figure}
\end{example}

\begin{example} Another tree with 54 vertices.
    \begin{figure}[H]
    \centering
    \includegraphics[width=\linewidth]{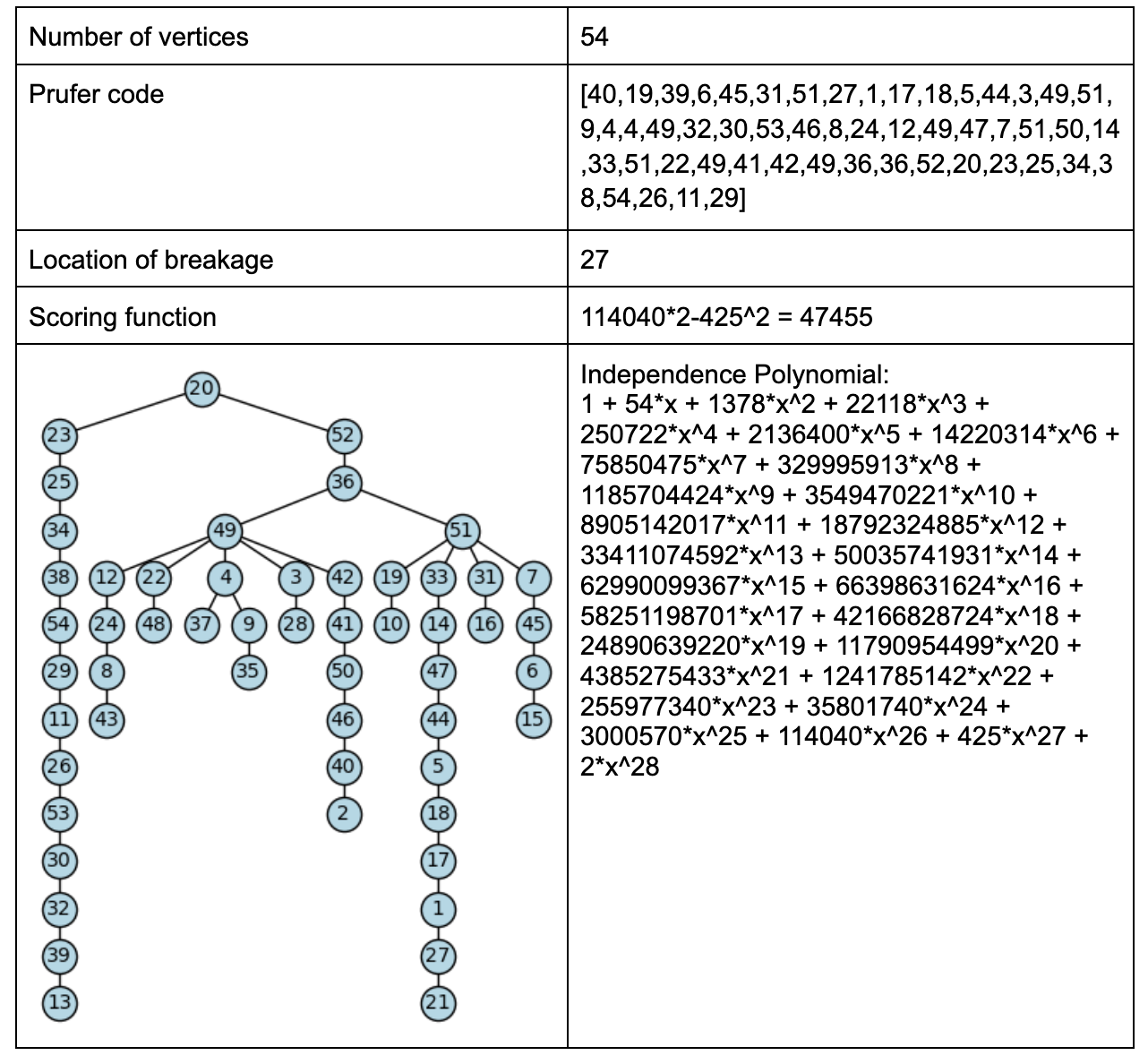}
    \end{figure}
\end{example}

\begin{example}\label{ex:T44} A tree with 37 vertices. This tree originally appeared in \cite{galvin2025trees} as $T_{4,4}$
    \begin{figure}[H]
    \centering
    \includegraphics[width=\linewidth]{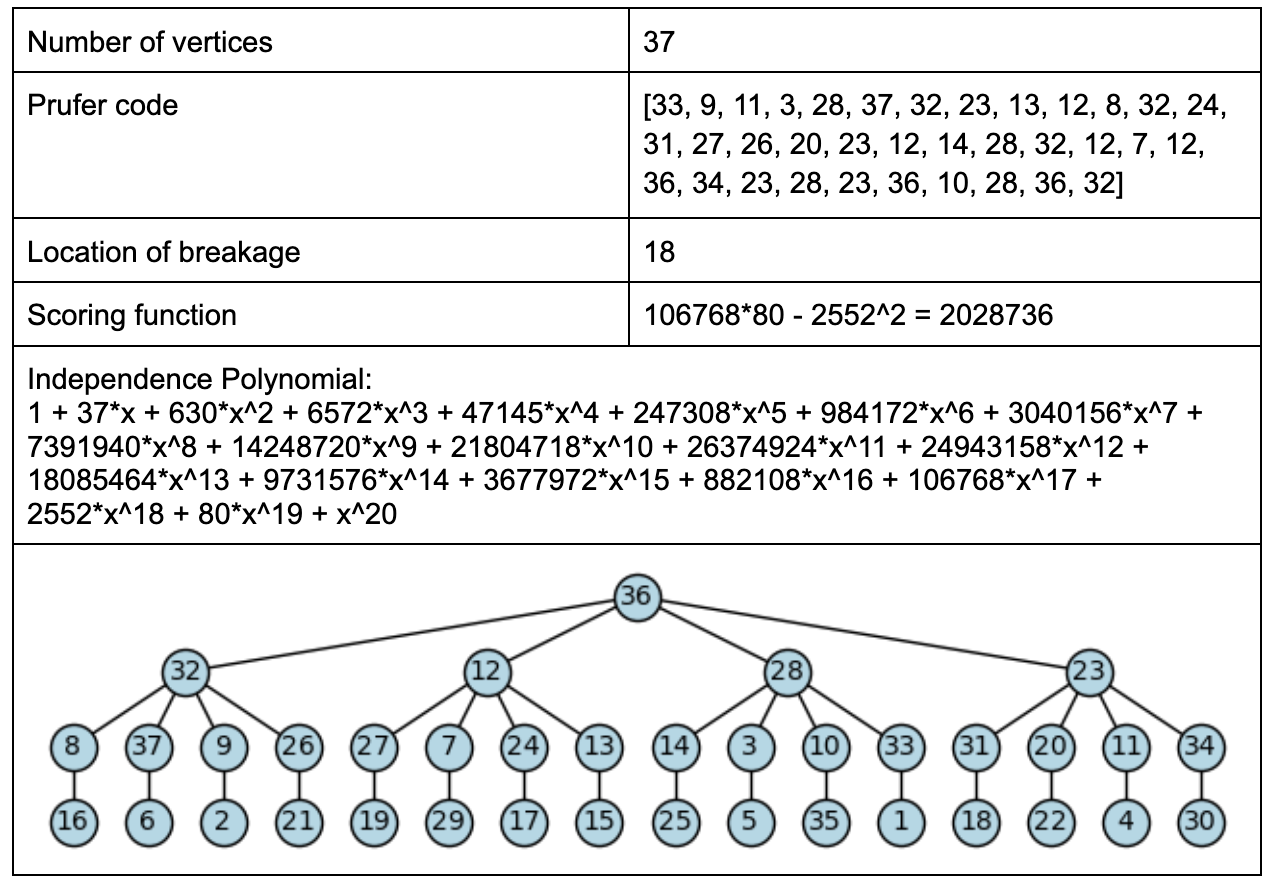}
    \end{figure}
\end{example}

\end{document}